 \documentclass[draft, 12pt, onecolumn]{IEEEtran}
\IEEEoverridecommandlockouts                              % This command is only
\title{\LARGE \bf
Information Relaxation and Dual Formulation of \\Controlled Markov Diffusions
}

\author{Fan~Ye and Enlu~Zhou  % <-this % stops a space

\thanks{F.~Ye and E.~Zhou are with the H. Milton Stewart School of Industrial \& Systems Engineering, Georgia Institute of Technology, Atlanta, GA, 30332, USA E-mail:{fye8}@gatech.edu \& {enlu.zhou}@isye.gatech.edu}
}

\usepackage[all]{xy}                    %支持交换图
\usepackage{amscd}
\usepackage{amsmath}
\usepackage{amssymb}                   % AMSLaTeX宏包
\usepackage{amsfonts}                  % AMSLaTeX宏包
\usepackage[amsmath,thmmarks]{ntheorem}% AMSLaTeX宏包.

\usepackage{graphicx,epsfig}                  %支持插图
\usepackage{mathrsfs}

 \usepackage{bbm}
\usepackage{latexsym}
\usepackage{lscape}
\usepackage{rotating}
\usepackage{graphicx}
\usepackage{mathptmx}

\renewcommand{\baselinestretch}{1.26}
\newtheorem{theorem}{Theorem}

\newtheorem{definition}{Definition}

\theoremstyle{remark} 
 \newtheorem{Proposition}{Proposition}
 \newtheorem{Remark}{Remark}
 \newtheorem{Assumption}{Assumption}

\begin{document}

\maketitle
%\IEEEpeerreviewmaketitle
%\thispagestyle{empty}
%\pagestyle{empty}
%
\pagestyle{headings}
\pagenumbering{arabic}
%%%%%%%%%%%%%%%%%%%%%%%%%%%%%%%%%%%%%%%%%%%%%%%%%%%%%%%%%%%%%%%%%%%%%%%%%%%%%%%%
\begin{abstract}
Information relaxation and duality in Markov decision processes  have been studied recently by several researchers with the goal to derive dual bounds on the value function. In this paper we extend this dual formulation to controlled Markov diffusions: in a similar way we  relax the constraint that the decision should be made based on the current information and impose a penalty to punish the access to the information in advance. We establish the weak duality, strong duality and complementary slackness results in a parallel way as those in Markov decision processes. We further explore the structure of the optimal penalties  and expose the connection between the optimal penalties for Markov decision processes and controlled Markov diffusions. We demonstrate the use of this dual representation   in a classic dynamic portfolio choice problem through a new class of penalties, which require little extra computation and produce small duality gap on the optimal value.
\end{abstract}
%\KEYWORDS{controlled diffusions, dynamic programming, information relaxation, duality, portfolio choice}

%%%%%%%%%%%%%%%%%%%%%%%%%%%%%%%%%%%%%%%%%%%%%%%%%%%%%%%%%%%%%%%%%%%%%%%%%%%%%%%%
\section{Introduction}\label{Section:Intro}
Markov decision processes (MDPs) and controlled Markov diffusions play a central role respectively in modeling discrete-time and continuous-time dynamic decision making problems under uncertainty, and hence  have wide applications in diverse fields such as engineering, operations research and economics. MDPs  and controlled Markov diffusions can be solved, in principle,  via dynamic programming and Hamlton-Jacobi-Bellman (HJB) equation, respectively. However, the exact computation of dynamic programming  suffers from the ``curse of dimensionality''- the size of the state space increases exponentially with the dimension of the state. Many  approximate dynamic programming  methods have been proposed for solving MDPs  to combat this curse of dimensionality, such as \cite{bertsekas:2007, chang:2007a, powell:2011, de2003linear}.  The HJB equation also rarely allows a closed-form solution, especially when the state space is of high dimension or there are constraints imposed on the controls. Several numerical methods have been developed including  \cite{kushner2001,Han:2011}; another standard numerical approach  is to discretize the time space, which reduces the original continuous-time problem to an MDP and hence the techniques of approximate dynamic programming can be applied.

It is worth noting that the approximate dynamic programming methods for solving MDPs often generate sub-optimal policies, and simulation under a sub-optimal policy leads to a lower bound (or upper bound) on the optimal expected reward (or cost). Though the accuracy of a sub-optimal policy is generally unknown, the lack of performance guarantee on a sub-optimal policy can be potentially addressed by  providing a dual bound, i.e., an upper bound (or lower bound) on the optimal expected reward (or cost). Valid and tight dual bounds based on a dual representation of MDPs  were recently developed by  \cite{rogers:2007} and \cite{brown:2010}.   The main idea of this duality approach is to relax the non-anticipativity constraints on decisions  but impose a penalty for getting access to the information in advance. In addition, this duality approach only involves pathwise deterministic optimization and therefore is well suited to Monte Carlo simulation, making it useful to evaluate the quality of sub-optimal policies in complex dynamic systems.

The dual formulation of MDPs is attractive in both theoretical and practical aspects. On one hand,  the idea of relaxing the non-anticipativity constraint on the control policies in  MDPs dates back  to at least  \cite{davis1995new}, as exposed by \cite{Haugh:2012}. In addition, the optimal penalty is not unique: for general problems we have the value function-based penalty developed by \cite{rogers:2007} and \cite{brown:2010}; for problems with convex structure there is an alternative optimal penalty, that is,  the gradient-based penalty, as pointed out by \cite{brown:2011}. On the other hand, in order to derive tight dual bounds, various approximation schemes based on  different optimal penalties have been proposed including \cite{brown:2010, brown:2011,Desai:2011, Ye2012parameterized}. We notice that this dual approach has found increasing applications in different fields, such as \cite{lai:2010, brown:2011,devalkar2011essays, moallemi2012dynamic,haugh2012dynamic}.

%especially when the dimension of the state space is high or there are constraints  on the control space,
The goal of this paper is to extend the information relaxation-based dual representation of MDPs to controlled Markov diffusions. Particularly, we intend to answer the following questions.
\begin{itemize}
\item Can we establish a similar framework of dual formulation for controlled Markov diffusions based on information relaxation as that for MDPs?
\item If the answer is yes, what is the form of the optimal penalty in the setting of controlled Markov diffusions? %Is the optimal penalty unique?
\item If certain optimal penalty exists, does its structure imply any computational advantage in deriving  dual bounds on the optimal value of practical problems?
%\item Is it possible to extend this dual representation to  more general continuous-time controlled Markov process?
\end{itemize}

%Since the quality of a numerical solution is hard to justify in many problems,  we are interested in deriving  a tight dual bound on the value function of a controlled Markov diffusion by formulating its dual representation. Around this topic some central questions are

The answer to the first question is yes, at least for a wide class of controlled Markov diffusions. To fully answer all the questions we present the information relaxation-based dual formulation of controlled Markov diffusions based on the technical machinery ``anticipating stochastic calculus''  (see, e.g., \cite{Ocone1989,nualart2006malliavin}).  We establish the weak duality, strong duality and complementary slackness results in a parallel way as those in the dual formulation of MDPs.  We investigate one type of  optimal penalties, i.e., the so-called ``value function-based  penalty'', to answer  the second question.  One key feature of the value function-based optimal penalty is that it can be written compactly as an Ito stochastic integral under the natural filtration generated by the Brownian motions. This compact expression potentially enables us to design sub-optimal penalties in simple forms and also facilitates the computation of the dual bound.  Then we emphasize on the computational aspect using the value function-based optimal penalty so as to answer the third question. A direct application  is illustrated by a classic dynamic portfolio choice problem with predictable returns and intermediate consumptions: we consider the numerical solution to a discrete-time model that is  discretized from a continuous-time model;  an effective class of penalties  that are easy to compute is proposed to derive dual bounds on the optimal value of  the discrete-time model.

%we propose a class of penalties based on the time discretization of the optimal value function-based penalties of the continuous-time model to compute dual bounds on the value function for the discrete-time model. These penalties make the dual approach much easier in terms of computation compared with the penalties directly derived from the discrete-time model.

It turns out that \cite{davis1992deterministic, davis1991anticipative, davis1989anticipative} have pioneered a series of related work for controlled Markov diffusions. They also adopted the approach of relaxing the future information and penalizing. In particular,  \cite{davis1992deterministic} proposed a Lagrangian approach for penalization, where the Lagrangian term plays essentially the same role as a penalty in our dual framework; in addition,  this Lagrangian term has a similar flavor as the gradient-based penalty proposed by \cite{brown:2011} for MDPs. The main difference of their work from ours is that we propose a more general framework that may incorporate their Lagrangian approach as a special case; the optimal penalty we develop in this paper is value function-based, which differs from  their proposed Lagrangian approach. In addition, their work is purely theoretical and does not suggest any computational method. In contrast,  we provide an example to demonstrate the practical use of the value function-based penalty.

Another closely-related literature focuses on the dual representation of the American option pricing problem (that is essentially  an optimal stopping problem) \cite{rogers:2002,haugh:2004,andersen:2004}.  In particular, the structure of the optimal martingale (i.e., the optimal penalty) under the diffusion process is investigated by \cite{Belomestny:2009,wang2010fast}, which leads to practical algorithms for fast computation of tight  upper bounds on the American option prices. The form of the optimal martingale also reflects its inherent relationship with the value function-based optimal penalty in the controlled diffusion setting.

We summarize our contributions as follows:
\begin{itemize}
\item We establish a dual representation of controlled Markov diffusions based on  information relaxation. We also explore the structure of the optimal penalty  and expose the connection between MDPs and controlled Markov diffusions. %, which can also incorporate the earlier work.
\item Based on the result of the dual representation of  controlled Markov diffusions, we demonstrate its practical use in a dynamic portfolio choice problem. In our numerical experiments the upper bounds on the optimal value show that  our proposed penalties are near optimal, comparing with the lower bounds induced by sub-optimal policies for the same problem.
%\item We propose the idea of parameterized penalties that can be used to tighten the dual bound.
\end{itemize}

%Finally, we note that our approximation scheme of the optimal penalty in the dynamic portfolio choice problem has a similar flavor of  using the martingale representation theorem to derive a near-optimal martingale in the continuous-time American option pricing (i.e., optimal stopping) problem \cite{Belomestny:2009} .

The rest of the paper is organized as follows. In Section \ref{Section:Dual_Representation}, we review the dual formulation of MDPs and derive the dual formulation of controlled Markov diffusions. In Section \ref{Section:Application}, we illustrate the dual approach and carry out  numerical studies in  a dynamic portfolio choice problem. Finally, we conclude with  future directions in Section \ref{Section:Conclusion}. We put some of the proofs and discussion of the connection between \cite{Belomestny:2009,wang2010fast} and our work in Appendix.
%, where we apply the dual approach to evaluate the upper bounds on the optimal expected utility

\section{Controlled Markov  Diffusions and Its Dual Representation} \label{Section:Dual_Representation}
We begin with a brief review of the dual framework on Markov Decision Processes that was first developed by \cite{rogers:2007} and \cite{brown:2010}. We then state the basic setup of the controlled Markov diffusion and its associated Hamilton-Jacobi-Bellman equation in Section \ref{Section:Controlled_Diffusion}. We develop the dual representation of controlled Markov diffusions and present the main results in Section \ref{Section:Dual_Diffusion}.% Some notations are abused across these subsections if they play a same role in different settings.

\subsection{Review of Dual Formulation of Markov Decision Processes}\label{Section:MDPs_formulation}
%We first briefly review the dual formulation of Markov decision processes(MDPs) based on information relaxation proposed in  \cite{brown:2010}.
Consider a finite-horizon MDP on the probability space $(\Omega,\mathcal{G},\mathbb{P})$. Time is indexed by $\mathcal{K}=\{0,1,\cdots,K\}$. Suppose $\mathcal{X}$ is the state space  and  $\mathcal{A}$ is the control space. The state $\{x_{k}\}$ follows the equation
\begin{align}
x_{k+1}&=f(x_{k},a_{k},v_{k+1}),~ k=0,1,\cdots,K-1, \label{state}
\end{align}
where $a_{k} \in \mathcal{A}_{k}$ is the control whose value is decided at time $k$, and  $\{v_{1},\cdots,v_{K}\}$ are independent  random variables for noise taking values in the set $\mathcal{V}$ with  known distributions.    The natural filtration is described by $\mathbb{G}=\{\mathcal{G}_{0},\cdots,\mathcal{G}_{K}\}$ with $\mathcal{G}_k\triangleq\sigma\{x_{0},v_{1}\cdots,v_{k}\}$; in particular, $\mathcal{G}=\mathcal{G}_{K}$.

Denote by $\mathbb{A}$ the set of all control strategies $\textbf{a}\triangleq (a_{0},\cdots,a_{K-1})$, i.e.,  each $a_k$ takes value in $\mathcal{A}$. Let $\mathbb{A}_{\mathbb{G}}$ be the set of control strategies that are adapted to the filtration $\mathbb{G}$, i.e., each $a_{k}$ is  $\mathcal{G}_{k}$-adapted. We also call $\textbf{a}\in\mathbb{A}_{\mathbb{G}}$ a  \emph{non-anticipative} policy. Given an $x_{0}\in \mathcal{X}$, the objective is to maximize the expected sum of intermediate rewards $\{g_k\}_{k=0}^{K-1}$ and final reward $\Lambda$   by selecting a non-anticipative policy $\textbf{a}\in\mathbb{A}_{\mathbb{G}}$:
\begin{align}
V_{0}(x_{0})&=\sup_{\textbf{a}\in \mathbb{A}_{\mathbb{G}}}J_{0}(x_{0};\textbf{a}), \notag \\
\text{where}~~J_{0}(x_{0};\textbf{a})&\triangleq \mathbb{E}\left[\sum_{k=0}^{K-1}g_{k}(x_{k},a_{k})+\Lambda(x_{K})\bigg|x_{0}\right]. \label{value_function}
\end{align}
%\begin{equation}
%V^{*}_{0}(x_{0})=\sup_{\textbf{a}\in \mathbb{A}_{\mathbb{F}}} \mathbb{E}[\sum_{t=0}^{T-1}g_{t}(x_{t},a_{t})+g_{T}(x_{T})|x_{0}], \label{value_function}
%\end{equation}
The expectation in (\ref{value_function}) is taken with respect to the random sequence $\textbf{v}=(v_{1},\cdots,v_{K})$.
%To avoid technical issues we assume that $V_{0}(x_{0},\textbf{a})$  has a uniform bound   for all $\textbf{a}\in\mathbb{A}_{\mathbb{F}}$.
The value function $V_{0}$ is a solution to the following  dynamic programming recursion:
\begin{align}
&V_{K}(x_{K})\triangleq \Lambda(x_{K}); \notag\\
&V_{k}(x_{k})\triangleq\sup_{a_{k}\in \mathcal{A}}\{g_{k}(x_{k},a_{k})+\mathbb{E}[V_{k+1}(x_{k+1})|x_{k},a_{k}]\},~k=K-1,\cdots,0.  \notag %\label{DP_recursion}
\end{align}

 Next we describe the dual formulation of the value function $V_{0}(x_{0})$. Here we only consider the \emph{perfect information relaxation}, i.e., we have full knowledge of the future randomness, since this relaxation is usually more applicable in practice.

Define $\mathbb{E}_{k,x}[\cdot]\triangleq \mathbb{E}[\cdot|x_k=x]$. Let $\mathcal{M}_{\mathbb{G}}(0)$ denote the set of \emph{dual feasible penalties} $M(\textbf{a},\textbf{v})$, which  do not penalize non-anticipative policies in expectation, i.e.,
 $$\mathbb{E}_{0,x}[M(\textbf{a},\textbf{v})]\leq 0 \text{~~for all~} x\in \mathcal{X}~ \text{and}~ \textbf{a}\in \mathbb{A}_{\mathbb{G}}.$$
  Denote by $\mathcal{D}$ the set of  real-valued functions on $\mathcal{X}$. Then we define an operator $\mathcal{L}:\mathcal{M}_{\mathbb{G}}(0)\rightarrow \mathcal{D}$:
\begin{align}
&\big(\mathcal{L}M\big)(x)=\mathbb{E}_{0,x}\left[\sup_{\textbf{a}\in \mathbb{A}} \left\{\sum_{k=0}^{K-1}g_{k}(x_{k},a_{k})+\Lambda(x_{K})-M(\textbf{a},\textbf{v})\right\}\right]. \label{dual_formulation}
\end{align}
Note that the supremum in (\ref{dual_formulation}) is over the set  $\mathbb{A}$ not the set  $\mathbb{A}_{\mathbb{G}}$, i.e.,  the control $a_k$ can be based on the future information.  The optimization problem inside the expectation in (\ref{dual_formulation}) is usually referred to as the \emph{inner optimization problem}. In particular, the right hand side of (\ref{dual_formulation}) is well suited to Monte Carlo simulation: we can simulate a realization of $\textbf{v}=\{v_{1},\cdots,v_{K}\}$ and solve the following inner optimization problem:

%\begin{align}
%I(x,M,\textbf{v})\triangleq\max_{\textbf{a},\textbf{x}}~~ &\sum_{k=0}^{K-1}g_{k}(x_{k},a_{k})+\Lambda(x_{K})- M(\textbf{a},\textbf{v})  \label{inner}\\
%\text{s.t.}~& x_{0}=x, \notag \\
%~ &x_{k+1}=f(x_{k},a_{k},v_{k+1}),~k=0,\cdots,K-1, \label{constraint}\\
%&a_{k}\in \mathcal{A}_{k},~k=0,\cdots,K-1,    \label{constraint2}
%\end{align}

\begin{subequations}\label{inner_opt}
\begin{align}
I(x,M,\textbf{v})\triangleq\max_{\textbf{a}}~~ &\sum_{k=0}^{K-1}g_{k}(x_{k},a_{k})+\Lambda(x_{K})- M(\textbf{a},\textbf{v})  \label{inner_opt1}\\
\text{s.t.}~& x_{0}=x, \notag \\
~ &x_{k+1}=f(x_{k},a_{k},v_{k+1}),~k=0,\cdots,K-1, \label{inner_opt1}\\
&a_{k}\in \mathcal{A}_{k},~k=0,\cdots,K-1,    \label{inner_opt3}
\end{align}
\end{subequations}
 which is in fact a \emph{deterministic} dynamic program. The optimal value  $I(x,M,\textbf{v})$ is an unbiased estimator of $(\mathcal{L}M)(x)$.

Theorem\ref{Theorem1:Duality_MDP} below establishes a strong  duality  in the sense that for all $x_0\in \mathcal{X},$
$$\sup_{\textbf{a}\in \mathbb{A}_{\mathbb{G}}}J_{0}(x_0;\textbf{a})= \inf_{M\in \mathcal{M}_\mathbb{G}(0)}\big(\mathcal{L}M\big)(x_0).$$
 In particular,
Theorem\ref{Theorem1:Duality_MDP}(a) suggests that $\mathcal{L}M(x_0)$ can be used to derive an upper bound on the value function $V_{0}(x_0)$  given any  $M\in \mathcal{M}_{\mathbb{G}}(0)$, i.e., $I(x_0,M,\textbf{v})$ is a high-biased estimator of $V_{0}(x_0)$ for all $x_0\in \mathcal{X}$; Theorem\ref{Theorem1:Duality_MDP}(b) states that the duality gap vanishes if the dual problem is solved by choosing $M$ in the form of (\ref{opt_penalty_MDPs}).
\begin{theorem}[Theorem 2.1 in \cite{brown:2010}]\label{Theorem1:Duality_MDP}
 \mbox{}
\begin{enumerate}
\item[(a)] (Weak Duality) For all  $M\in\mathcal{M}_{\mathbb{G}}(0)$  and all $x\in \mathcal{X},$~ $V_{0}(x)\leq (\mathcal{L}M)(x).$
\item[(b)] (Strong Duality) For all $x\in \mathcal{X}$, $V_{0}(x)= (\mathcal{L}M^{*})(x)$, where \\
\begin{equation}
M^{*}(\textbf{a},\textbf{v})=\sum_{k=0}^{K-1}\left(V_{k+1}(x_{k+1})-\mathbb{E}[V_{k+1}(x_{k+1})|x_{k},a_{k}]\right). \label{opt_penalty_MDPs}
\end{equation}
\end{enumerate}
\end{theorem}
\begin{Remark}
 \mbox{}
\begin{enumerate}
\item Note that the right hand side of (\ref{opt_penalty_MDPs}) is a function of (\textbf{a},\textbf{v}), since $\{x_k\}$ depend on $(\textbf{a},\textbf{v})$ through the equation (\ref{state}).
%\item The reason that $M\in \mathcal{M}_{\mathbb{G}}(0)$ is called a (dual feasible) penalty function becomes clear: if the relaxation of the requirement on the non-anticipative policies is penalized  by using a proper function in $\mathcal{M}_{\mathbb{G}}(0)$, then  the value function $V_{0}$ can be recovered via the dual approach due to the strong duality result.

\item Note that the optimal penalty $M^{*}(\textbf{a},\textbf{v})$ is  the sum of a  $\mathbb{G}$-martingale difference sequence when $\textbf{a}\in \mathbb{A}_{\mathbb{G}}$; therefore, $M^{*}(\textbf{a},\textbf{v})\in \mathcal{M}_{\mathbb{G}}(0).$ Since  $M^{*}$ depends on the value function $\{V_k\}$, it is referred to as the \emph{value function-based penalty}.
%. So we call an $M\in\mathcal{M}_{\mathbb{G}}(0)$ a dual feasible penalty.
\end{enumerate}
\end{Remark}

The optimal penalty (\ref{opt_penalty_MDPs}) that  achieves the strong duality involves the value function $\{V_k\}$, and hence is intractable in practical problems. In order to obtain tight dual bounds, a natural idea is to derive sub-optimal penalty functions based on  a good approximate value function $\{\hat{V}_k\}$  or some sub-optimal policy $\hat{\textbf{a}}$.  Methods based on these ideas have been successfully implemented in the American option pricing problems by \cite{rogers:2002,haugh:2004,andersen:2004}, and also  in \cite{brown:2010,lai:2010,devalkar2011essays}.

\subsection{Controlled Markov Diffusions and Hamilton-Jacobi-Bellman Equation}\label{Section:Controlled_Diffusion}

This subsection is concerned with the control of Markov diffusion processes. Applying the Bellman's principle of dynamic programming  leads to a second-order nonlinear partial differential equation, which is referred to as the Hamilton-Jacobi-Bellman equation. For a comprehensive treatment on this topic we refer the readers to \cite{Fleming:2006}.

Let us consider a $\mathbb{R}^n$-valued controlled Markov diffusion process $(x_{t})_{0\leq t\leq T}$ driven by an $m$-dimensional Brownian motion $(w_t)_{0\leq t\leq T}$ on a probability space $(\Omega,\mathcal{F},\mathbb{P})$, following the stochastic differential equation (SDE):
\begin{equation}
dx_{t}=b(t,x_t,u_t)dt+\sigma(t,x_t)dw_t,~~ 0 \leq t \leq T \label{state_SDE},
\end{equation}
where the control $u_t$ takes value in $\mathcal{U}\subset \mathbb{R}^{d_u}$ ($d_u\in \mathbb{N}$), while $b$ and $\sigma$ are functions $b:[0,T]\times \mathbb{R}^{n}\times \mathcal{U}\rightarrow \mathbb{R}^{n}$ and $\sigma:[0,T]\times \mathbb{R}^{n}\rightarrow \mathbb{R}^{n\times m}$. The natural (augmented) filtration generated by the Brownian motions is denoted by $\mathbb{F}=\{\mathcal{F}_{t}, 0 \leq t \leq T\}$ with $\mathcal{F}=\mathcal{F}_{T}$. In the following $\parallel\cdot \parallel$ denotes the Euclidean norm.

\begin{definition}
A control strategy $\textbf{u}=(u_s)_{s\in[t, T]}$ is called an admissable strategy at time $t$ if
\begin{enumerate}
\item $\textbf{u}=(u_s)_{s\in[t, T]}$ is an $\mathbb{F}$-progressively measurable process  taking values in $\mathcal{U}$ (i.e., $\textbf{u}$ is a non-anticipative policy), and satisfying $\mathbb{E}[\int_{t}^{T}||u_s||^2ds]<\infty$;
%\item The equation (\ref{state_SDE}) has a pathwise unique  solution $(x_s)_{s\in[t,T]}$ that is  $\mathbb{F}$-progressively measurable and has continuous sample paths given $x_t=x\in \mathbb{R}^{n}$;
\item  $\mathbb{E}_{t,x}[\sup_{s\in[t,T]}||x_s||^2]<\infty$, where $\mathbb{E}_{t,x}[\cdot]\triangleq\mathbb{E}[\cdot|x_t=x]$.
\end{enumerate}
The set of  admissible strategies  at time $t$ is denoted by  $\mathcal{U}_{\mathbb{F}}(t)$.
\end{definition}
With the following standard technical conditions imposed on $b$ and $\sigma$, the SDE (\ref{state_SDE}) admits a unique pathwise solution when  $\textbf{u}\in \mathcal{U}_{\mathbb{F}}(0)$, i.e., $(x_t)_{0\leq t \leq T}$ is $\mathbb{F}$-progressively measurable and has continuous sample paths almost surely given $x_0=x\in \mathbb{R}^{n}$.
\begin{Assumption}
$b$ and $\sigma$ are continuous on their domains, respectively,  and for some constants $C_{1}, C_{2},$ and $C_{\sigma}>0$,
\begin{enumerate}
\item $\parallel b(t,x,u)\parallel+ \parallel \sigma(t,x)\parallel \leq C_{1}(1+\parallel x \parallel+\parallel u\parallel)$  for all $(t,x)\in \bar{Q}$ and $u\in\mathcal{U}$;

\item $\parallel b(t,x,u)-b(s,y,u)\parallel + \parallel \sigma(t,x)-\sigma(s,y)\parallel \leq C_2(|t-s|+\parallel x-y\parallel)$ for all  $(t,x),(s,y)\in \bar{Q}$  and $u\in\mathcal{U}$.

\item $\xi^{\top}(\sigma\sigma^{\top})(t,x) \xi\geq C_{\sigma}\parallel \xi\parallel^2$ for all $(t,x)\in [0,T]\times Q$ and $\xi\in \mathbb{R}^{n}$.
\end{enumerate}
\end{Assumption}
Let $Q=[0,T)\times \mathbb{R}^{n}$ and  $\bar{Q}=[0,T]\times \mathbb{R}^{n}$. We define the functions  $\Lambda:\mathbb{R}^{n}\rightarrow \mathbb{R}$  and $g:\bar{Q}\times \mathcal{U}\rightarrow \mathbb{R}$ as the final reward and intermediate reward, respectively. Assume that  $\Lambda$ and $g$ satisfy the following polynomial growth conditions.
\vspace{-2mm}
\begin{Assumption}
For some constants $C_{\Lambda},c_{\Lambda},C_{g},c_{g}>0$,
  \begin{enumerate}
 \item $|\Lambda(x)|\leq C_{\Lambda}\left(1+\parallel x\parallel^{c_{\Lambda}}\right)$  for all $x\in \mathbb{R}^n $;
 \item $|g(t,x,u)|\leq C_{g}\left(1+\parallel x\parallel^{c_{g}}+\parallel u\parallel^{c_{g}}\right)$ for all $(t,x)\in \bar{Q}$. %for some constants $C_{g},,c_{g}>0$
  \end{enumerate}
  \end{Assumption}
Then we introduce the reward functional
\begin{equation*}
J(t,x;\textbf{u})\triangleq \mathbb{E}_{t,x}\left[\Lambda(x_{T})+\int_{t}^{T}g(s,x_s,u_s)ds\right].
\end{equation*}
 Given an initial condition $(t,x)\in Q$, the objective is to maximize  $J(t,x,u)$  over all the controls $\textbf{u}$ in $\mathcal{U}_{\mathbb{F}}(t)$:
%; currently we  simply assume they are continuous and bounded functions on their domains.
\begin{equation}
V(t,x)=\sup_{\textbf{u}\in \mathcal{U}_{\mathbb{F}}(t)} J(t,x;\textbf{u}). \label{value_function_SDE}
\end{equation}
Here we abuse the notations of the state $x$, the rewards $\Lambda$ and $g$, and the value function $V$, since they play the same roles as those in MDPs.

Let $C^{1,2}(Q)$ denote the space of function $L(t,x):Q\rightarrow \mathbb{R}$ that is continuously differentiable in (i.e., $C^{1}$) in $t$ and twice continuously differentiable (i.e., $C^2$) in $x$ on $Q$.  For $L\in C^{1,2}(Q)$, define a partial differential operator $A^{u}$ by
\begin{align*}
A^{u}L(t,x)\triangleq &L_t(t,x)+L^{\top}_x(t,x)b(t,x,u)+\frac{1}{2}\text{tr}\left(L_{xx}(t,x)\left(\sigma\sigma^{\top}\right)(t,x)\right),
\end{align*}
where $L_t$, $L_x$, and $L_{xx}$ denote the $t$-partial derivative, the gradient and the Hessian with respect to $x$ respectively, and  $\left(\sigma\sigma^{\top}\right)(t,x)\triangleq\sigma(t,x)\sigma^{\top}(t,x)$. Let $C_{p}(\bar{Q})$ denote the set of function $L(t,x):\bar{Q}\rightarrow \mathbb{R}$ that is continuous on  $\bar{Q}$ and satisfies a polynomial growth condition in $x$, i.e.,
$$|L(t,x)|\leq C_{L}(1+\parallel x\parallel^{c_{L}})$$
for some constants $C_{L}$ and $c_{L}$. The following well-known verification theorem provides a sufficient condition for the value function and an optimal control strategy using Bellman's principle of dynamic programming. %(see, e.g.,Theorem 3.1 in \cite{Fleming:2006}).
%Suppose $V\in C^{1,2}([0,T)\times \mathbb{R}^{n})$. Then $V(t,x)$  necessarily should satisfy the Hamilton-Jacobi-Bellman equation:

\begin{theorem}[Verification Theorem, Theorem 4.3.1 in \cite{Fleming:2006}]\label{Theorem2:Verification}
Suppose Assumptions 1 and 2 hold, and $\bar{V}\in  C^{1,2}(Q)\cap C_{p}(\bar{Q})$ satisfies
\begin{equation}
\sup_{u\in \mathcal{U} }\{g(t,x,u)+A^{u}\bar{V}(t,x)\}=0  ~\text{for}~(t,x)\in Q, \label{HJB}
\end{equation}
and $\bar{V}(T,x)=\Lambda(x).$  Then\\
(a) $J(t,x;\textbf{u})\leq \bar{V}(t,x)$ for any $\textbf{u}\in \mathcal{U}_{\mathbb{F}}(t)$ and any $(t,x)\in \bar{Q}$.\\
(b) If there exists a function $u^*:\bar{Q}\rightarrow \mathcal{U}$ such that
\begin{equation}
g(t,x,u^{*}(t,x))+A^{u^{*}(t,x)}\bar{V}(t,x)=\max_{u\in \mathcal{U} }\{g(t,x,u)+A^{u}\bar{V}(t,x)\}=0 \label{HJB1}
\end{equation}
 for all $(t,x)\in Q$ and if the control strategy defined as $\textbf{u}^{*}=(u^{*}_t)_{t\in[0,T]}$ with $u^{*}_t\triangleq u^*(t,x_t)$ is admissible at time $0$ (i.e., $\textbf{u}^{*}\in \mathcal{U}_{\mathbb{F}}(0)$), then
 \begin{enumerate}
\item $\bar{V}(t,x)=V(t,x)=\sup_{\textbf{u}\in \mathcal{U}_{\mathbb{F}}(t)} J(t,x;\textbf{u}).$ for all $(t,x)\in \bar{Q}$.
\item $\textbf{u}^{*}$ is an optimal control strategy, $i.e.$, $V(0,x)=J(0,x;\textbf{u}^{*})$.
\end{enumerate}
 \end{theorem}
Equation (\ref{HJB}) is the well-known HJB equation associated with the problem (\ref{state_SDE})-(\ref{value_function_SDE}). %However, the existence of $\bar{V}\in C^{1,2}(Q)$ in Theorem \ref{Theorem2:Verification} requires many technical assumptions that might not be true in practice. For example, the HJB equation is usually assumed to be of  \emph{uniformly parabolic} type if there exists $c_{\sigma}>0$ such that for all $(t,x,u)\in Q\times\mathcal{U}$ and $\xi\in \mathbb{R}^{n}$,
%$$\xi^{\top}(\sigma\sigma^{\top})(t,x,u) \xi\geq c_{\sigma}\parallel \xi\parallel^2.$$
%Otherwise, a classic solution $\bar{V}\in C^{1,2}(Q)$ may not be expected and we need to interpret the value function as a viscosity solution to the HJB equation (see, e.g., \cite{Fleming:2006}). % The other technical conditions to ensure the existence of $\bar{V}\in C^{1,2}(Q)$ can be found in IV.4 of \cite{Fleming:2006}.

 %For ease of understanding we provide a proof in Appendix \ref{Section:AppendixVerification}.

\subsection{Dual Representation of Controlled Markov Diffusions}\label{Section:Dual_Diffusion}
In this subsection we present the information relaxation-based dual formulation of controlled Markov diffusions. In a similar way we relax the constraint that the decision at every time instant should be made based on the past information and impose a penalty to punish the access to future information. We will establish the weak duality, strong duality and complementary slackness results for controlled Markov diffusions, which parallel the results in MDPs.  The value function-based optimal penalty is also characterized to motivate the practical use of our dual formulation, which will be demonstrated in Section \ref{Section:Application}.

We consider the perfect information relaxation, i.e.,  we  can foresee all the future randomness generated by the Brownian motion so that the decision made at any time $t\in [0,T]$ is  based on the information set $\mathcal{F}=\mathcal{F}_T$.  To expand the set of the feasible controls, we use  $\mathcal{U}(t)$  to denote the set of measurable $\mathcal{U}$-valued control strategies at time $t$, i.e., $\textbf{u}=(u_{s})_{s\in[t,T]}\in \mathcal{U}(t)$ if $\textbf{u}$ is $\mathcal{B}([t,T])\times\mathcal{F}$-measurable and   $u_s$ takes value in $\mathcal{U}$ for $s\in[t,T] $, where $\mathcal{B}([t,T])$ is the Borel $\sigma$-algebra on $[t,T]$. In particular, $\mathcal{U}(0)$ can be viewed as the counterpart of $\mathbb{A}$ introduced in Section \ref{Section:MDPs_formulation} for MDPs.

Unlike the case of MDPs, the first technical problem we have to face  is to define a solution of (\ref{state_SDE}) with an anticipative control $\textbf{u}\in  \mathcal{U}(0)$.  Since it involves the concept of ``anticipating stochastic calculus'' and Stratonovich integral, we postpone the technical details to Appendix \ref{Section:Appendix:anticipating}, where  we use the decomposition technique to define the solution of an anticipating  SDE  following \cite{davis1992deterministic}, \cite{Ocone1989}.

%This is also the reason why we restrict the dependence of $\sigma(t,x,u)$ in (\ref{state_SDE}) only on $t$ and $x$.

 Right now we assume that given a control strategy $\textbf{u}\in \mathcal{U}(0)$ there exists a unique solution $(x_t)_{t\in[0,T]}$ to (\ref{state_SDE})  that is $\mathcal{B}([0,T])\times\mathcal{F}$-measurable. Next we consider the set of penalty functions in the setting of controlled Markov diffusions.   Suppose  $h(\textbf{u},\textbf{w})$ is a function depending on a control strategy $\textbf{u}\in\mathcal{U}(0)$ and a sample path of Brownian motion  $\textbf{w}\triangleq(w_t)_{t\in[0,T]}$. We define the set $\mathcal{M}_{\mathbb{F}}(0)$  of \emph{dual feasible penalties} $h(\textbf{u},\textbf{w})$ that  do not penalize non-anticipative policies in expectation, i.e.,
 % We use $\textbf{w}=(w_t)_{t\in[0,T]}$ to denote a sample path of Brownian motion.
\begin{equation*}
\mathbb{E}_{0,x}[h(\textbf{u},\textbf{w})]\leq 0 \text{~~for all~} x\in \mathbb{R}^{n} \text{~and~} \textbf{u}\in\mathcal{ U}_{\mathbb{F}}(0).
\end{equation*}
%such that $(x_s)_{s\in[t,T]}$ has a unique strong solution.
In the following we will show  $\mathcal{M}_{\mathbb{F}}(0)$ parallels the role of $\mathcal{M}_{\mathbb{G}}(0)$ for MDPs  in the dual formulation of controlled Markov diffusions.

%plays the same role in the dual formulation of controlled Markov diffusions as $\mathcal{M}_{\mathbb{G}}(0)$ does in the dual formulation of  MDPs.

With an arbitrary choice of $h\in \mathcal{M}_{\mathbb{F}}(0)$, we can determine an upper bound on (\ref{value_function_SDE}) with $t=0$  by relaxing the constraint on the adaptiveness of control strategies.

\begin{Proposition}[Weak Duality] \label{Prop:weak_duality}
If $h\in \mathcal{M}_{\mathbb{F}}(0)$, then for all  $x\in \mathbb{R}^{n},$
\begin{eqnarray}
\sup_{\textbf{u}\in \mathcal{U}_{\mathbb{F}}(0)} J(0,x;\textbf{u}) \leq \mathbb{E}_{0,x}\left[\sup_{\textbf{u}\in \mathcal{U}(0)}\left\{\Lambda(x_{T})+\int_{0}^{T}g(t,x_t,u_t)dt-h(\textbf{u},\textbf{w})\right\}\right]. \label{weak_duality}
\end{eqnarray}
\end{Proposition}

\begin{IEEEproof}
For any $\bar{\textbf{u}}\in \mathcal{U}_{\mathbb{F}}(0)$,
\begin{align*}
J(0,x;\bar{\textbf{u}})=&\mathbb{E}_{0,x}\left[\Lambda(x_{T})+\int_{0}^{T}g(t,x_t,\bar{u}_t)dt\right]\\
\leq&  \mathbb{E}_{0,x}\left[\Lambda(x_{T})+\int_{0}^{T}g(t,x_t,\bar{u}_t)dt-h(\bar{\textbf{u}},\textbf{w})\right]\\
\leq& \mathbb{E}_{0,x}\left[\sup_{\textbf{u}\in \mathcal{U}(0)}\left\{\Lambda(x_{T})+\int_{0}^{T}g(t,x_t,u_t)dt-h(\textbf{u},\textbf{w})\right\}\right].
\end{align*}
Then inequality (\ref{weak_duality}) can be obtained by taking the supremum over $\bar{\textbf{u}}\in \mathcal{U}_{\mathbb{F}}(0)$ on the left hand side of the last inequality.
\end{IEEEproof}

The optimization problem inside the conditional expectation in (\ref{weak_duality}) is the counterpart of  (\ref{inner_opt})  in the context of controlled Markov diffusions: an entire path of $\textbf{w}$ is known beforehand (i.e., \emph{perfect information relaxation}),  and the objective function depends on a specific trajectory of $\textbf{w}$. Therefore, it is a deterministic and path-dependent optimal control problem parameterized by \textbf{w}. We also call it an \emph{inner optimization problem}, and  the expectation term on the right hand side of (\ref{weak_duality}) is a \emph{dual bound} on the value function $V(0,x)$.  References \cite{davis1992deterministic, davis1989anticipative, davis1991anticipative} have conducted a series of  research on this problem  under the name of ``anticipative stochastic control''. In particular, one of the special cases they have considered is $h=0$, which means the future information is accessed  without any penalty;  \cite{davis1992deterministic} characterized the value of the perfect information relaxation.  We would expect that  the dual bound associated with the zero penalty can be very loose as that in MDPs.  The evaluation of the dual bound is well suited to Monte Carlo simulation: we can generate a sample path of $\textbf{w}$ and solve the inner optimization problem in (\ref{weak_duality}), the solution of which is a high-biased estimator of $V(0,x)$.

An interesting case is when we choose
\begin{equation}
h^{*}(\textbf{u},\textbf{w})=\Lambda(x_{T})+\int_{0}^{T}g(t,x_t,u_t)dt-V(0,x).\label{optimal_penalty_1}
\end{equation}
Note that $h^{*}\in\mathcal{M}_{\mathbb{F}}(0)$, since by the definition of $V(0,x),$

$$\mathbb{E}_{0,x}\left[\Lambda (x_{T})+\int_{0}^{T}g(s,x_s,u_s)ds\right]\leq V(0,x) \text{~~for all~} x\in \mathbb{R}^{n} \text{~and~} \textbf{u}\in \mathcal{U}_{\mathbb{F}}(0).$$

We also note that by plugging $h=h^{*}$ in the inner optimization problem in (\ref{weak_duality}), the objective value of which is independent of $\textbf{u}$ and it is always equal to $V(0,x)$. So the following strong duality result is obtained.

\begin{theorem}[Strong Duality]\label{Theorem3:StrongDuality}
For all $x\in \mathbb{R}^{n}$,
\begin{equation}
\sup_{\textbf{u}\in \mathcal{U}_{\mathbb{F}}(0)} J(0,x;\textbf{u})=\inf_{h\in \mathcal{M}_{\mathbb{F}}(0)}\left\{\mathbb{E}_{0,x}\left[\sup_{\textbf{u}\in \mathcal{U}(0)}\left\{\Lambda(x_{T})+\int_{0}^{T}g(t,x_t,u_t)dt-h(\textbf{u},\textbf{w})\right\}\right]\right\}.\label{strong_duality_1}
\end{equation}
The minimum of the right hand side of (\ref{strong_duality_1}) can always be achieved by choosing an $h\in \mathcal{M}_{\mathbb{F}}(0)$ in the form of (\ref{optimal_penalty_1}).
\end{theorem}

\begin{IEEEproof}
 According to the weak duality, the left side of (\ref{strong_duality_1}) should be less than or equal to the right side of (\ref{strong_duality_1}); the equality is achieved by choosing $h=h^{*}$ in (\ref{optimal_penalty_1}).
\end{IEEEproof}

Due to the strong duality result, the left side of (\ref{strong_duality_1}) is referred to as the \emph{primal problem} and the right  side of (\ref{strong_duality_1}) is referred to as the \emph{dual problem}.   If $\textbf{u}^{\star}$ is a control strategy that achieves the supremum  in the primal problem, and  $h^{\star}$ is a dual feasible penalty that achieves the infimum in the dual problem, then they are optimal solutions to the primal and dual problems, respectively. The ``complementary slackness condition''  in the next theorem, which parallels the result in the discrete-time problem (Theorem 2.2 in \cite{brown:2010}), characterizes such a pair $(\textbf{u}^{\star}, h^{\star})$.

\begin{theorem}[Complementary Slackness]\label{Theorem4:Complementary}
Given $\textbf{u}^{\star}\in \mathcal{U}_{\mathbb{F}}(0)$ and $h^{\star}\in\mathcal{M}_{\mathbb{F}}(0)$, a sufficient and necessary condition for $\textbf{u}^{\star}$ and $h^{\star}$ being optimal to the primal and dual problem respectively is that
$$\mathbb{E}_{0,x}[h^{\star}(\textbf{u}^{\star},\textbf{w})]=0,$$
and
\begin{align}
&\mathbb{E}_{0,x}\left[\Lambda(x^{\star}_{T})+\int_{t}^{T}g(s,x^{\star}_s,u^{\star}_s)ds-h^{\star}(\textbf{u}^{\star},\textbf{w})\right] \notag \\
=&\mathbb{E}_{0,x}\left[\sup_{\textbf{u}\in \mathcal{U}(0)}\left\{\Lambda(x_{T})+\int_{0}^{T}g(s,x_s,u_s)ds-h^{\star}(\textbf{u},\textbf{w})\right\}\right],\label{ComplementarySlackness}
\end{align}
where $x^{\star}_t$ is the solution of (\ref{state_SDE}) using the control strategy $\textbf{u}^{\star}=(u^{\star}_t)_{t\in[0,T]}$ on $[0,t)$ with the initial condition $x^{\star}_0=x$.
\end{theorem}
\begin{IEEEproof}
We first consider sufficiency. Let $\textbf{u}^{\star}\in \mathcal{U}_{\mathbb{F}}(0)$ and $h^{\star}\in \mathcal{M}_{\mathbb{F}}(0)$. We assume $\mathbb{E}_{0,x}[h^{\star}(\textbf{u}^{\star},\textbf{w})]=0$ and (\ref{ComplementarySlackness}) holds. Then by the weak duality, $\textbf{u}^{\star}$ and $h^{\star}$ should be optimal to the primal and dual problem, respectively.

Next we consider necessity.  Let $\textbf{u}^{\star}\in \mathcal{U}_{\mathbb{F}}(0)$ and $h^{\star}\in \mathcal{M}_{\mathbb{F}}(0)$.  Then we have
\begin{align*}
&\mathbb{E}_{0,x}\left[\sup_{\textbf{u}\in \mathcal{U}(0)}\left\{\Lambda(x_{T})+\int_{0}^{T}g(t,x_t,u_t)dt-h^{\star}(\textbf{u},\textbf{w})\right\}\right] \notag\\
\geq &\mathbb{E}_{0,x}\left[\Lambda(x^{\star}_{T})+\int_{t}^{T}g(t,x^{\star}_t,u^{\star}_t)dt-h^{\star}(\textbf{u}^{\star},\textbf{w})\right]\\
\geq &J(0,x;\textbf{u}^{\star}).
\end{align*}
The last inequality holds due to $h^{\star}\in \mathcal{M}_{\mathbb{F}}(0)$. Since we know  $\textbf{u}^{\star}$ and $h^{\star}$ are optimal to the primal and dual problem respectively,  then by the strong duality result
$$J(0,x;\textbf{u}^{\star})= \mathbb{E}_{0,x}\left[\sup_{\textbf{u}\in \mathcal{U}(0)}\left\{\Lambda(x_{T})+\int_{0}^{T}g(t,x_t,u_t)dt-h^{\star}(\textbf{u},\textbf{w})\right\}\right],$$
which implies all the inequalities above are equalities. Therefore, we know $\mathbb{E}_{0,x}[h^{\star}(\textbf{u}^{\star},\textbf{w})]=0$ and (\ref{ComplementarySlackness}) holds.
\end{IEEEproof}

Here we have the same interpretation on complementary slackness condition as that in the dual formulation of  MDPs: if the penalty  is optimal to the dual problem, the decision maker will be satisfied with an optimal non-anticipative control strategy even if she is able to choose any anticipative control strategy. Clearly, if an optimal control strategy $\textbf{u}^{*}$ to the primal problem (\ref{state_SDE})-(\ref{value_function_SDE}) does exist (see, e.g., Theorem \ref{Theorem2:Verification}(b)), then $\textbf{u}^{*}$ and $h^{*}(\textbf{u},\textbf{w})$ defined in (\ref{optimal_penalty_1}) is a pair of the optimal solutions to the primal and dual problems. However, we note that the optimal penalty in the form of (\ref{optimal_penalty_1}) is intractable as it depends on the exact value of $V(0,x)$. The next theorem characterizes the form of another optimal penalty, which motivates the numerical approximation scheme that will be illustrated in Section \ref{Section:Application}.

 %We fully develop the relevant results in Theorem\ref{Theorem5:StrongDual_IdealPenalty} in Appendix \ref{Section:Appendix:value_penalty} to avoid obscuring the main message and intuition here.

%The optimal penalty  $h^{*}(\textbf{u},\textbf{w})$ can be written compactly as an Ito stochastic integral,  when it is evaluated at $\textbf{u}=\textbf{u}^*$. A natural question would be  whether $\int_{0}^{T}V_x(t,x_t)^{\top}\sigma(t,x_t)dw_t $ plays the role of an optimal penalty  in (\ref{strong_duality_1})  as $M^{*}(\textbf{a},\textbf{v})$ does in Theorem \ref{Theorem1:Duality_MDP} achieving the strong duality. Unfortunately, $\int_{0}^{T}V_x(t,x_t)^{\top}\sigma(t,x_t)dw_t $ is not even a well-defined object in terms of an Ito stochastic integral, when $\textbf{u}$  is not adapted to $\mathbb{F}$. To fix this problem, we also need the machinery of ``anticipating stochastic calculus''. However, we can still provide a concise answer here, that is, there exists an alternative optimal penalty that coincides with  $\int_{0}^{T}V_x(t,x_t)^{\top}\sigma(t,x_t)dw_t$ when  $\textbf{u}\in \mathcal{U}_{\mathbb{F}}(0)$. We fully develop the relevant results in Theorem \ref{Theorem5:StrongDual_IdealPenalty} in Appendix \ref{Section:Appendix:value_penalty}.  In the following proposition  we formalize one of the main results in Theorem \ref{Theorem5:StrongDual_IdealPenalty}, which also guides the numerical approximation scheme  that will be illustrated in Section \ref{Section:Application}.

\begin{theorem}[Value Function-Based Penalty]\label{Theorem5:StrongDual_IdealPenalty}
Suppose that the value function $V(t,x)$ for the problem (\ref{state_SDE})-(\ref{value_function_SDE})  satisfies  the assumptions in Theorem\ref{Theorem2:Verification}(b), and $\textbf{y}=(t,x_t)_{t\in[0,T]}$ satisfies the conditions in Proposition \ref{Theorem7:Ito_Stratonovich} in
 Appendix \ref{Section:Appendix:anticipating} (i.e., the Ito formula for Stratonovich integral (\ref{Stratonovich_ito_formula}) is valid for $F=V(t,x)$ and $\textbf{y}=(t,x_t)_{t\in[0,T]}$), where $(x_t)_{t\in[0,T]}$ is the solution to (\ref{state_SDE}) with  $\textbf{u}\in \mathcal{U}(0)$.  For $\textbf{u}\in \mathcal{U}(0)$, define
\begin{align}
 &h^{*}_v(\textbf{u},\textbf{w})\triangleq  \sum_{i=1}^m\int_0^T\left[V^{\top}_x(t,x_t)\sigma^i(t,x_t)\right]\circ dw_t^i \notag \\
&-\frac{1}{2}\int_0^T\left[V_x^{\top}(t,x_t)\left(\sum_{i=1}^{m}\sigma^i_x\sigma^i(t,x_t)\right)+\text{tr}\left(V_{xx}(t,x_t)(\sigma\sigma^{\top})(t,x_t)\right)\right]\,dt.  \label{optimal_penalty_3}
\end{align}
Then
\begin{enumerate}
\item If $\textbf{u}\in \mathcal{U}_{\mathbb{F}}(0)$, (\ref{optimal_penalty_3}) reduces to the  form
\begin{align}
h^{*}_v(\textbf{u},\textbf{w})=\int_0^TV^{\top}_x(t,x_t)\sigma(t,x_t)\,dw_t,\label{optimal_penalty_2}
\end{align}
and $h^{*}_v(\textbf{u},\textbf{w})\in  \mathcal{M}_{\mathbb{F}}(0)$.
\item The strong duality holds in
\begin{align*}
V(0,x)=\mathbb{E}_{0,x}\left[\sup_{\textbf{u}\in \mathcal{U}(0)}\left\{\Lambda(x_{T})+\int_{0}^{T}g(t,x_t,u_t)dt-h^{*}_v(\textbf{u},\textbf{w})\right\}\right]. %\label{strong_duality}
\end{align*}

Moreover, the following equalities hold almost surely with $x_0=x$
\begin{align}
V(0,x)=&\sup_{\textbf{u}\in \mathcal{U}(0)}\left\{\Lambda(x_{T})+\int_{0}^{T}g(t,x_t,u_t)dt-h^{*}_v(\textbf{u},\textbf{w})\right\}\label{almost_sure}\\
=&\Lambda(x_{T}^{*})+\int_{0}^{T}g(t,x^{*}_t,u^{*}_t)dt-h^{*}_v(\textbf{u}^{*},\textbf{w}),\label{almost_sure_2}
\end{align}
where $(x^{*}_t)_{t\in[0,T]}$ is the solution of  (\ref{state_SDE}) using the optimal control $\textbf{u}^{*}=(u^{*}_t)_{t\in[0,T]}$ (defined in Theorem\ref{Theorem2:Verification}(b)) on $[0,t)$ with the initial condition $x^{*}_0=x$.
\end{enumerate}
\end{theorem}

%\begin{Proposition}\label{Prop:optimal_penalty_1}
%Suppose  the value function $V(t,x)$ in (\ref{value_function_SDE}) satisfies all the assumptions in Theorem\ref{Theorem2:Verification}(b). Then under the technical condition specified in Theorem\ref{Theorem5:StrongDual_IdealPenalty}, there is an optimal solution to the dual problem, i.e., an optimal penalty $h^{*}_v(\textbf{u},\textbf{w})\in\mathcal{M}_{\mathbb{F}}(0)$  in the form of
%\begin{equation}
%h^{*}_v(\textbf{u},\textbf{w})=\int_{0}^{T}V^{\top}_x(t,x_t)\sigma(t,x_t)dw_t~~\text{for}~ \textbf{u}\in \mathcal{U}_{\mathbb{F}}(0), \label{optimal_penalty_2}
%\end{equation}
% where $x_t$ is the solution of (\ref{state_SDE}) using the control $\textbf{u}=(u_t)_{t\in[0,T]}$  on $[0,t)$ with the initial condition $x_0=x$.
%\end{Proposition}

Since the value functions $\{V(t,x), 0\leq t\leq T\}$ are unknown in real applications, (\ref{optimal_penalty_2})  implies that if an approximate value function $\{\hat{V}(t,x), 0\leq t\leq T\}$ is differentiable with respect to $x$, then heuristically,  $h^{*}_v$  can be approximated by   $\hat{h}_v(\textbf{u},\textbf{w})\triangleq \int_{0}^{T}\hat{V}^{\top}_x(t,x_t)\sigma(t,x_t)dw_t$  at least for $\textbf{u}\in \mathcal{U}_{\mathbb{F}}(0)$.  Noting that $\{\int_0^t\hat{V}^{\top}_x(s,x_s)\sigma(s,x_s)dw_s\}_{0\leq t\leq T}$ is an $\mathbb{F}$-martingale if $\textbf{u}\in \mathcal{U}_{\mathbb{F}}(0)$ (assuming that $\hat{V}^{\top}_x(t,x)\sigma(t,x)$ satisfies the polynomial growth condition in $x$); therefore,
$\mathbb{E}_{0,x}[\hat{h}_v(\textbf{u},\textbf{w})]=0 \text{~~for all~} x\in \mathbb{R}^{n} \text{~and~} \textbf{u}\in\mathcal{ U}_{\mathbb{F}}(0).$ As a result, $\hat{h}_v(\textbf{u},\textbf{w})\in\mathcal{M}_{\mathbb{F}}(0)$, i.e., $\hat{h}$ is dual feasible, which means that $\hat{h}_v$ can be used to derive an upper bound on the value function $V(0,x)$ through (\ref{weak_duality}). Hence, in terms of the approximation scheme implied by the form of the optimal penalty,  Theorem\ref{Theorem5:StrongDual_IdealPenalty} presents a \emph{value function-based penalty} that  can be viewed as the continuous-time analogue of $M^{*}(\textbf{a},\textbf{v})$ in (\ref{opt_penalty_MDPs}).

It is revealed by  the complementary slackness condition in both discrete-time (Theorem 2.2 in \cite{brown:2010}) and continuous-time (Theorem\ref{Theorem4:Complementary}) cases that any optimal penalty has zero expectation evaluating at  an optimal  policy;  as a stronger version, the value function-based optimal penalty in both cases  assign zero expectation to all non-anticipative polices  (note that $M^{*}$ in (\ref{opt_penalty_MDPs}) is a sum of martingale differences  under the original filtration $\mathbb{G}$).
%Under the condition that $\hat{h}_v(\textbf{u},\textbf{w})$ is close to $h^*_v(\textbf{u},\textbf{w})$, i.e., $E[\hat{h}_v(\textbf{u},\textbf{w})-h^*_v(\textbf{u},\textbf{w})]$ for all $\textbf{u}\in \mathcal{U}(0)$,  then the upper bounds induced by $\hat{h}_v(\textbf{u},\textbf{w})$ could be tight.

 Intuitively, we can interpret the  strong duality achieved by the value function-based penalty as to offset the path-dependent randomness in the inner optimization problem;  then the optimal control to the inner optimization problem   coincides with that to  the original stochastic control problem in the expectation sense, which is reflected by the proof of Theorem\ref{Theorem5:StrongDual_IdealPenalty} in Appendix \ref{Section:Appendix:value_penalty} for controlled Markov diffusions. In Appendix \ref{Section:Appendix:optimal_stopping} we briefly review the dual representation of the optimal stopping problem, where an analogous result of  Theorem\ref{Theorem5:StrongDual_IdealPenalty} exists  provided the evolution of the state is modelled as a diffusion process.

\section{Dynamic Portfolio Choice Problem} \label{Section:Application}

We illustrate the practical use of the dual formulation of controlled Markov diffusions, especially the value function-based optimal penalty developed in Theorem\ref{Theorem5:StrongDual_IdealPenalty}, in a classic  dynamic portfolio choice problem with predictable returns and intermediate consumptions (see, e.g., \cite{samuelson1969lifetime,merton1969lifetime,merton1971optimum}).  Since most portfolio choice problems of practical interest cannot be solved analytically, various numerical methods have been developed  including the martingale approach \cite{cvitanic:2003,he:1991},  state-space discretization methods \cite{tauchen:1991,balduzzi:1999}, and approximate dynamic programming methods \cite{brandt:2005,Han:2011}. These methods all produce sub-optimal policies, and it is  not difficult to obtain  lower bounds on the optimal expected utility by Monte Carlo simulation under these policies; on the other hand, an upper bound is constructed by \cite{haugh:2006} and  \cite{brown:2011}  respectively based on the work by \cite{cvitanic:1992} and  \cite{brown:2010}. The gap between the lower bound and the upper bound can be used to justify the performance of a candidate policy.

%Though some methods bear the property of asymptotic convergence, the accuracy of numerical solutions with limited computational power  cannot be measured.
%constructed an upper bound on the  optimal value  based on the dual formulation of the constrained portfolio choice problem proposed by \cite{cvitanic:1992} and the information relaxation duality method proposed by \cite{brown:2010,rogers:2007}.

In this section we solve a \emph{discrete-time} dynamic portfolio choice problem  that is discretized from  a  continuous-time model (see, e.g., \cite{cvitanic:1992,liu:2007}). We consider the time-discretization as it is a common approach to numerically solve the continuous-time problem, and  the decisions of investment only occur at discrete-time points.   We focus on generating upper bounds on  the optimal expected utility of the discrete-time problem using the information relaxation dual approach. In particular, we propose a new class of penalties for the discrete-time problem by discretizing  the value function-based optimal penalties of the continuous-time problem. These penalties make the inner optimization problem much easier to solve  compared with the  penalties  that directly approximates the optimal penalty of the discrete-time model. We demonstrate the  effectiveness of our method in computing dual bounds through numerical experiments.

%The information relaxation technique based on the value function-based penalty employed in this section can be naturally generalized to other controlled diffusion problems. However, the main challenge often lies in solving the inner optimization problem with both effectiveness and efficiency, which is a central issue in the practical use of this dual approach. Our dynamic portfolio example is one of the first few examples of applying the information relaxation-based dual approach to continuous-state MDPs.

%So far there is very few examples directly wokring towards, The dynamic portfolio choice example in this section serves as the

%
%we propose a new approach to generate penalty functions that avoid evaluating any conditional expectation and solve an inner optimization problem for each random sequence.  We demonstrate the  effectiveness of our method in computing dual bounds through numerical experiments with a risk-free asset and three risky assets.

\subsection{The Portfolio Choice Model}\label{section:sub:model}
We first consider a continuous-time financial market with finite horizon $[0,T]$, which is built on the probability space $(\Omega, \mathcal{F}, \mathbb{P})$. There are one risk-free asset and $n$ risky assets that the investor can invest on. The prices of the risk-free asset and risky assets are denoted by $S_{t}^{0}$ and  $S_{t}=(S_{t}^1,\cdots,S_{t}^n)^{\top}$, respectively, and the instantaneous asset returns depend on the $m$-dimensional state variable $\phi_t$:
\begin{align}
dS_{t}^{0}&=r_{f}S_{t}^{0}dt \notag\\
dS_t&= S_t\bullet(\mu_tdt+\sigma_t dz_t), \label{continuous_asset} \\
d\phi_t&=\mu^{\phi}_tdt+\sigma^{\phi,1}_tdz_t+\sigma^{\phi,2}_td\tilde{z}_t, \label{continuous_market_state}
\end{align}
where $r_f$  is the instantaneous risk-free rate of return, and $\textbf{z}\triangleq(z_{t})_{0\leq t\leq T}$ and $\tilde{\textbf{z}}\triangleq(\tilde{z}_{t})_{0\leq t\leq T}$  are two independent  standard Brownian motions  that are of dimension $n$ and $d$, respectively; the  drift vector $\mu_t=\mu(t,\phi_t)$  and the diffusion matrix $\sigma_t=\sigma(t,\phi_t)$   in (\ref{continuous_asset}) are of dimension $n$ and $n\times n$, where the symbol  $\bullet$  denotes the component-wise multiplication of two vectors; the terms $\mu^{\phi}_t=\mu^{\phi}(t,\phi_t)$, $\sigma^{\phi,1}_t=\sigma^{\phi,1}(t,\phi_t)$, $\sigma^{\phi,2}_t=\sigma^{\phi,2}(t,\phi_t)$ in (\ref{continuous_market_state})  are of dimension $m$, $m\times n$, and $m \times d$, respectively.

We denote the filtration by $\mathbb{F}=\{\mathcal{F}_{t},0\leq t \leq T\}$, where $\mathcal{F}_t$ is generated by the Brownian motions $\{(z_{s},\tilde{z}_{s}), 0\leq s \leq t \}$.

%We will use $\Sigma_t$ to denote  $\sigma_t\sigma^{\top}_t$,  the covariance matrix of the instantaneous return.

%$$dS_{t}^{0}=r_{f}S_{t}^{0}dt.$$
%The price vector of $n$ risky assets is  denoted by  $S_{t}=(S_{t}^1,\cdots,S_{t}^n)^{\top}$  and it follows a geometric Brownian motion
%\begin{align}
%dS_t=(\mu_t\circ S_t) dt+\sigma_t^{\top}S_t dz_t,\label{continuous_asset}
%\end{align}
%where $\textbf{z}\triangleq(z_{t})_{0\leq t\leq T}$ is an $n$-dimensional standard Brownian motion. The  drift vector $\mu_t=\mu(t,\phi_t)$  and the diffusion matrix $\sigma_t=\sigma(t,\phi_t)$   are of dimension $n$ and $n\times n$, respectively,  where  $\phi_t$ is an $m$-dimensional market state variable that follows another diffusion process
%\begin{align}
%d\phi_t=\mu^{\phi}_tdt+\sigma^{\phi,1}_tdz_t+\sigma^{\phi,2}_td\tilde{z}_t,\label{continuous_market_state}
%\end{align}
%where $\mu^{\phi}_t=\mu^{\phi}(t,\phi_t)$, $\sigma^{\phi,1}_t=\sigma^{\phi,1}(t,\phi_t)$, $\sigma^{\phi,2}_t=\sigma^{\phi,2}(t,\phi_t)$ are of dimension $m$, $m\times n$, and $m \times d$, respectively, and  $\tilde{\textbf{z}}\triangleq(\tilde{z}_{t})_{0\leq t\leq T}$ is another $d$-dimensional standard Brownian motion independent of $\textbf{z}.$
%Denote the filtration by $\mathbb{F}=\{\mathcal{F}_{t},0\leq t \leq T\}$, where $\mathcal{F}_t$ is generated by $\{(z_{s},\tilde{z}_{s}), 0\leq s \leq t \}$. The covariance matrices $\sigma_t\sigma^{\top}_t$ is denoted by  $\Sigma_t$.

Let $\pi_t=(\pi^{1}_t,\cdots,\pi^{n}_t)^{\top}$ and $\tilde{c}_t$ denote the fraction of wealth  invested in  $n$ risky assets and  the instantaneous rate of  consumption, respectively. The total wealth $W_{t}$ of a portfolio that consists of  the $n$ risky assets and one risk-free asset evolves according to
\begin{align}
dW_{t}=&W_t\left[\pi_t^{\top}\left(\mu_t dt+\sigma_t dz_t\right)+r_f\left(1-\pi_{t}^{\top}\textbf{1}_{n}\right)dt-\tilde{c}_tdt\right] \notag\\
=&W_t\left(\pi_t^{\top}(\mu_t-r_f \textbf{1}_{n})+r_f  -\tilde{c}_t\right)dt+W_t\pi_t^{\top}\sigma_t dz_t, \label{continuous_wealth}
\end{align}
%=&\sum_{i=1}^{n}\pi^{i}_t (\mu^{i}_tdt+\sum_{j=1}^{n}\sigma^{ij}_tdz_{t}^{j})+r_f(1-\textbf{1}^{T}\pi_{t})dt,
where $\textbf{1}_{n}$ is the $n$-dimensional all-ones vector. The control process $\textbf{u}\triangleq(u_t)_{0\leq t \leq T}$ with $u_t\triangleq (\pi_t,\tilde{c}_t)$  is  an admissible strategy in the sense that
\begin{enumerate}
\item The control $\textbf{u}$ is $\mathbb{F}$-progressively measurable and $\mathbb{E}[\int_{0}^{T}||u_t||^2dt]<\infty$;
\item  $W_t> 0$, $\tilde{c}_t\geq0$, and $\int_0^{T}W_t\tilde{c}_tdt <\infty$ a.s.; %and $\pi_t$ satisfying a square integrability condition $$\int_{0}^{T}\parallel\pi_t\parallel^{2}dt<\infty;$$
\item $u_t\in \mathcal{U}$, where $\mathcal{U}$ is a closed convex set in $\mathbb{R}^{n+1}$.
\end{enumerate}
We still use  $\mathcal{U}_{\mathbb{F}}(t)$ to denote the set of admissible strategies at time $t$ and we will specify the control space $\mathcal{U}$ later. Suppose that $U$ is a strictly increasing and concave utility function (see, e.g., \cite{luenberger1997investment}). The investor's objective is to maximize the weighted sum of the expected utility of  the intermediate consumption and  the final wealth:
\begin{align}
V(t,\phi_t,W_t)=\sup_{\textbf{u}\in \mathcal{U}_{\mathbb{F}}(t)}\mathbb{E}&\bigg[\int_{t}^{T}\alpha\beta^{ s}U\left(\tilde{c}_sW_s\right)ds +(1-\alpha)\beta^{T}U(W_{T})\bigg|\phi_t, W_t\bigg], \label{Value_portfolio}
\end{align}
where $\beta\in[0,1)$ is the discount factor, and $\alpha\in[0,1]$ indicates the relative importance of the intermediate consumption.

%The standard stochastic control approach is to solve the HJB equation:
% \begin{align}
% 0&=\sup_{u\in \mathcal{U}}\big[\beta^{ t}U_1(\tilde{c}_tW_t)+A^u V(t,\phi_t,W_t)\big], \label{continous_HJB}
% \end{align}
%where
%\begin{align*}
%A^u V(\cdot)=&V_{t}(\cdot)+V^{\top}_{\phi}(\cdot)\mu^{\phi}(t,\phi_t)+V_{W}(\cdot) W_t(\pi_{t}^{\top}(\mu_t-r_f)+r_f-\tilde{c}_t)\\
%&+\frac{1}{2}\text{tr}[V_{\phi\phi}(\cdot)(\Sigma_1^{\phi}(t,\phi_t)+\Sigma_2^{\phi}(t,\phi_t))] \notag
% +V_{\phi W}^{\top}(\cdot)\sigma_{1}^{\phi}(t,\phi_t)(W_t\pi_t^{\top}\sigma_t)^{\top}\\
% &+\frac{1}{2}V_{WW}(\cdot)W^{2}_{t}\pi_t^{\top}\Sigma_{t}(t,\phi_t)\pi_t
%\end{align*}
%with the terminal  condition $V(T,\phi_{T}, W_{T})=(1-\alpha)\beta ^{T}U_2(W_{T})$.

The value function (\ref{Value_portfolio}) sometimes admits an analytic solution, for example, under the assumption that  $\mu_t$ is a constant vector  and $\sigma_t$ is a constant matrix in  (\ref{continuous_asset}), and there is no constraint on $u_t=(\pi_t,\tilde{c}_t)$. A recent progress on the analytic tractability of  (\ref{Value_portfolio}) can be found in \cite{liu:2007}.  However, (\ref{Value_portfolio}) usually does not have an analytic result when there is a position constraint on $\pi_t$.

Considering that the investment and consumption can only take place in a finite number of times in the real world,  we discretize the continuous-time problem (\ref{continuous_market_state})-(\ref{Value_portfolio}). Suppose the decision takes place at equally spaced times $\{0=t_0,t_1\cdots, t_{K}\}$ such that  $K=T/\delta$, where $\delta=t_{k+1}-t_k$ for $k=0,1,\cdots,K-1$. We simply denote the time grids by $\{0,1,\cdots,K\}$.   Note that  (\ref{continuous_asset}) is equivalent to
\begin{equation*}
d\log(S_t)=\left(\mu_t-\frac{1}{2} \cdot \text{Pdiag}\left(\Sigma_t\right)\right) dt+\sigma_tdz_t,
\end{equation*}
where $\text{Pdiag}(\Sigma_t)$ denotes an $n$-dimensional vector that is the principal diagonal of $\Sigma_t=\sigma_t\sigma^{\top}_t$,  the covariance matrix of the instantaneous return. That is to say, $S_{k+1}=R_{k+1}\bullet S_k$ with  distribution $\log(R_{k+1})\sim N(\int_{k\delta}^{(k+1)\delta}(\mu_s-\frac{1}{2}\sigma_s^2)ds, \int_{k\delta}^{(k+1)\delta} \Sigma_s ds).$ Hence, we can discretize (\ref{continuous_market_state}),(\ref{continuous_asset}), and (\ref{continuous_wealth}) as follows:
\begin{subequations}\label{discrt_model}
\begin{align}
\phi_{k+1}&=\phi_{k}+\mu^{\phi}_k\delta +\sigma^{\phi,1}_k \sqrt{\delta}Z_{k+1}+\sigma^{\phi,2}_k\sqrt{\delta}\tilde{Z}_{k+1}, \label{discrt_market_state}\\
\log(R_{k+1})&=\left(\mu_{k}-\frac{1}{2}\sigma_{k}^2\right)\delta  + \sigma_k\sqrt{\delta} Z_{k+1}, \label{discrt_R}\\
W_{k+1}&= W_k\left(R_{k+1}^{\top}\pi_k\right)+W_k\left(1-\textbf{1}_{n}^{\top}\pi_k\right)R_f-W_kc_k, \notag\\
&=W_k\left(R_f+(R_{k+1}-R_{f}\textbf{1}_{n})^{\top}\pi_k-c_k\right), \label{discrt_wealth}
\end{align}
\end{subequations}
where $\{(Z_{k},\tilde{Z}_{k}), k=1,\cdots,K\}$ is a sequence of identically and independently distributed standard Gaussian random vectors. In particular, we use $R_f\triangleq 1+r_f\delta$ and the decision variable $c_k$  to approximate  $e^{r_f\delta}$ and $\tilde{c}_k\delta$ due to the discretization procedure.

%$\mu^{\phi}_k=\mu^{\phi}(k,\phi_{k})$, $\sigma^{\phi,1}_k=\sigma^{\phi,1}(k,\phi_{k})$, $\sigma^{\phi,2}_k=\sigma^{\phi,2}(k,\phi_{k})$, $\mu_k=\mu(\phi_k)$, $\sigma_k=\sigma(\phi_k)$

Here we abuse the notations $\phi,W,$ and $\pi$ in the continuous-time and discrete-time settings. However, the subscripts make them easy  to distinguish: the subscript $t\in[0,T]$ is used in the continuous-time model, while $k=0,\cdots,K$ is used in the discrete-time model.

Denote the filtration of the process (\ref{discrt_model}) by $\mathbb{G}=\{\mathcal{G}_{0},\cdots,\mathcal{G}_{K}\}$, where $\mathcal{G}_k$ is generated by $\{(Z_{j},\tilde{Z}_{j}), j=0,\cdots,k\}$. In our numerical examples we assume that short sales and borrowing are not allowed, and the consumption cannot exceed the amount of the risky-free asset.  Then the constraint, on the control $a_k\triangleq(\pi_k,c_k)$ for the discrete-time problem,  can be defined as
\begin{equation}
 \mathcal{A}\triangleq\{(\pi,c)\in \mathbb{R}^{n+1}|\pi\geq0, c\geq 0, c \leq R_f(1-\textbf{1}_n^{\top}\pi)\}. \label{discrt_constraint}
\end{equation}
Since  $c_k$  is used to approximate $\tilde{c}_k\delta$, (\ref{discrt_constraint}) corresponds to a control set for the continuous-time model, which is  defined as
$$\mathcal{U}\triangleq \{(\pi,\tilde{c})\in \mathbb{R}^{n+1}|\pi\geq0, \tilde{c}\geq0, \tilde{c}\leq R_f(1-\textbf{1}_n^{\top}\pi)/\delta \}.$$

Let  $\mathbb{A}_{\mathbb{G}}$ again denote the set of $\mathcal{A}$-valued control strategies $\textbf{a}\triangleq (a_{1},\cdots,a_{K-1})$ that are adapted to the filtration $\mathbb{G}$. The discretization of (\ref{Value_portfolio}) serves as the value function to the discrete-time problem:
\begin{equation}
H_0(\phi_0,W_0)=\sup_{\textbf{a}\in \mathbb{A}_{\mathbb{G}}}\mathbb{E}_0\left[\sum_{k=0}^{K-1} \alpha\beta^{k\delta}U(c_k W_k)\delta+(1-\alpha)\beta^{K\delta}U(W_{K})\right], \label{discrt_value}
\end{equation}
which can be solved  via dynamic programming:
\begin{align}
H_{K}(\phi_K,W_K)&=(1-\alpha) \beta^{K\delta}U(W_K); \notag\\
H_{k}(\phi_k,W_k)&=\sup_{a_k\in\mathcal{A}}\left\{\alpha\beta^{k\delta}U(c_kW_k)\delta+ \mathbb{E}_k\left[H_{k+1}\left(\phi_{k+1},W_{k+1}\right)\right]\right\} . \label{dynamic_program}
\end{align}

We will focus on \emph{solving the discrete-time model }(\ref{discrt_model})-(\ref{discrt_value}), which is discretized from the continuous-time model (\ref{continuous_market_state})-(\ref{Value_portfolio}). Though our methods proposed later can be applied on general utility functions, for the purpose of illustration we consider the utility functions of the constant relative risk aversion (CRRA) type with coefficient $\gamma>0$, $i.e$, $U(x)=\frac{1}{1-\gamma}x^{1-\gamma}$, which are widely used in economics and finance. Since the utility functions are of CRRA type, both value functions (\ref{Value_portfolio}) and (\ref{discrt_value}) have simplified structures. To be specific, the value function to the continuous-time problem can be written as the factorization (see, e.g., \cite{liu:2007})
\begin{equation}
V(t,\phi_t,W_t)=\beta^{t}W_t^{1-\gamma}\tilde{J}(t,\phi_t), \label{value_function_SDE_structure}
\end{equation}
where $\tilde{J}(T,\phi_T)=(1-\alpha)/(1-\gamma),$ and
$$\tilde{J}(t, \phi)=\sup_{\textbf{u}\in \mathcal{U}_{\mathbb{F}}(t)}\mathbb{E}\left[\int_{t}^{T}\beta^{ s-t}\frac{\alpha}{1-\gamma}\left(\tilde{c}_sW_s\right)^{1-\gamma}ds+\beta^{ T-t}\frac{1-\alpha}{1-\gamma}W_{T}^{1-\gamma}\bigg|\phi_t=\phi,W_t=1\right];$$
and the value function to the discrete-time problem, due to the factorization scheme, can be written as
\begin{equation}
H_{k}(\phi_k,W_k)=\beta^{k\delta}W_k^{1-\gamma}J_{k}(\phi_k), \label{discrt_value_fact}
\end{equation}
 where $J_k$, the discrete-time reward functional,  is defined recursively as $J_K(\phi_K)=(1-\alpha)/(1-\gamma)$ and
\begin{equation}
J_k(\phi_k)=\sup_{(\pi_k,c_k)\in \mathcal{A}}\left\{\frac{\alpha}{1-\gamma}c_k^{1-\gamma}\delta+\beta^{\delta}\mathbb{E}\left[ \big(R_f+(R_{k+1}-R_{f})^{\top}\pi_k-c_k\big)^{1-\gamma}J_{k+1}(\phi_{k+1})|\phi_k \right]\right\} . \label{recursion_portfolio}
\end{equation}

It can be seen that the structure of the value functions to both continuous-time model and discrete-time model are similar: they
can be decomposed as a product of a function of the wealth  $W$ and a function of the market state variable $\phi$. If $\delta$ is small, $\tilde{J}(k\delta,\phi)$ and $J_{k}(\phi)$ may be close to each other. As a byproduct of this decomposition, another feature of the dynamic portfolio choice problem
with CRRA utility function  is that the optimal asset allocation and consumption $(\pi_{t},\tilde{c}_t)$ in continuous-time model are independent of the wealth $W_t$ given $\phi_t$ (respectively, the optimal $(\pi_{k},c_k)$ in discrete-time model are independent of the wealth $W_k$ given $\phi_k$). So the dimension of the
state space in (\ref{dynamic_program}) is actually the dimension of $\phi_k$. A number of numerical methods have been developed to solve the discrete-time model based on the recursion (\ref{recursion_portfolio}) including the state-space discretization approach \cite{tauchen:1991, balduzzi:1999}, and a simulation-based method \cite{brandt:2005}.

%We should also note from (\ref{factor}) that the dynamic portfolio choice at time $k$ is identical to the myopic (single-period)
%portfolio choice if $R_{k+1}$ and $\phi_{k+1}$ are independent conditional on $\phi_k$, since the conditional expectation term can be factorized as
%a product of two (conditional expectation) terms; the optimal $\pi_t$ is independent of the second term.  A typical example is when $\sigma^{\phi}_{1}=0$.

\subsection{Penalties and Dual Bounds}
In this subsection, we compute upper bounds on the optimal value $H_0$ of the discrete-time (and continuous-state) model (\ref{discrt_model})-(\ref{discrt_value}) based on the dual approach for MDPs in Theorem 1. We illustrate how to  generate two dual feasible penalties: one directly approximates the value function-based penalty of the discrete-time problem, while the other one is derived by discretizing the value function-based penalty of the continuous-time problem (\ref{continuous_market_state})-(\ref{Value_portfolio}). We discuss why the latter approach is more desirable to  generate upper bounds on  $H_0$  in terms of computational tractability of the inner optimization problem.

Throughout this subsection  we assume that an approximate function of $J_{k}(\phi)$, say $\hat{J}_{k}(\phi)$ (therefore, $\hat{H}_{k}(\phi_k,W_k)\triangleq W_k^{1-\gamma}\hat{J}_{k}(\phi_k)$ is an approximation of $H_k$), and an approximate policy $\hat{\textbf{a}}\in \mathbb{A}_{\mathbb{G}}$ are available.  We do not require that $\hat{\textbf{a}}$ should be derived from $\hat{J}_{k}(\phi)$ or vice versa; in other words, they can be obtained using different approaches.  We first describe the information relaxation dual approach of MDPs in the context of our portfolio choice problem. We focus on the perfect information relaxation that assumes the investor can foresee the future uncertainty $\textbf{Z}=(Z_1,\cdots,Z_{K})$ and $\tilde{\textbf{Z}}=(\tilde{Z}_1,\cdots,\tilde{Z}_{K})$, i.e., all the market states  and  returns of the risky assets.  A function $M(\textbf{a},\textbf{Z},\tilde{\textbf{Z}})$  is  a \emph{dual feasible} penalty in the setting of dynamic portfolio choice problem if for any $(\phi_0,W_0)$,
 \begin{equation}
 \mathbb{E}\left[M(\textbf{a},\textbf{Z},\tilde{\textbf{Z}})|\phi_0,W_0\right]\leq 0 \text{~~for all~~} \textbf{a}\in \mathbb{A}_{\mathbb{G}}. \label{portfolio_penalty_set}
 \end{equation}
 Let $\mathcal{M}_\mathbb{G}(0)$ denote the set of  all dual feasible penalties. For $M\in \mathcal{M}_\mathbb{G}(0)$ we  define  $\mathcal{L}M$ as a function of  $(\phi_0,W_{0})$:
\begin{align}
(\mathcal{L}M)(\phi_0,W_{0})=\mathbb{E}\bigg[\sup_{\textbf{a}\in \mathbb{A} } \{\sum_{k=0}^{K-1}\alpha \beta^{k\delta}U(c_k W_k)\delta&+(1-\alpha)\beta^{K\delta}U(W_{K}) -M(\textbf{a},\textbf{Z},\tilde{\textbf{Z}})\}\bigg|\phi_0,W_{0}\bigg]. \label{dual_portfolio}
\end{align}
Based on Theorem\ref{Theorem1:Duality_MDP}(a), $(\mathcal{L}M)(\phi_0,W_{0})$ is an upper bound on $H_0(\phi_0,W_{0})$ for any $M\in \mathcal{M}_\mathbb{G}(0)$ .% Note that the supremum in (\ref{dual_formulation}) is over the set of feasible strategies $\mathbb{A}$ not the set of  nonanticipative policies $\mathbb{A}_{\mathbb{F}}$.

%The optimization problem inside the expectation in (\ref{dual_formulation}) is referred to as the \emph{inner optimization problem} in the rest of this section.

 To ease the inner optimization problem, we introduce equivalent  decision variables $\Pi_k=W_k\pi_k$ and $C_k=W_k c_k$, which can be interchangeably used with $\pi_k$ and $c_k$.  We still use $\textbf{a}$  to denote an admissable strategy, though in terms of $(\Pi_k,C_k)$ now. Then we can rewrite the inner optimization problem inside the conditional expectation in (\ref{dual_portfolio}) as follows:
\begin{subequations}\label{inner_portfolio_opt}
\begin{align}
I(\phi_0,W_{0},M,\textbf{Z},\tilde{\textbf{Z}})&\triangleq \max_{\Pi,C,W} \left\{\sum_{k=0}^{K-1} \alpha\beta^{k\delta}U(C_k)\delta+(1-\alpha)\beta^{K\delta}U(W_{K})-M(\textbf{a},\textbf{Z},\tilde{\textbf{Z}})\right\} \label{inner_porfolio_obj}\\
\text{s.t.}~&\phi_{k+1}=\phi_{k}+\mu^{\phi}_k\delta +\sigma^{\phi,1}_k\sqrt{\delta}Z_{k+1}+\sigma^{\phi,2}_k\sqrt{\delta}\tilde{Z}_{k+1}, \label{inner_porfolio_constraint4}\\
&\log(R_{k+1})=(\mu_{k}-\frac{1}{2}\sigma_{k}^2)\delta  + \sigma_k\sqrt{\delta} Z_{k+1}, \label{inner_porfolio_constraint5} \\
&W_{k+1}=W_{k}R_f+(R_{k+1}-R_f\textbf{1}_n)^{\top}\Pi_k-C_k,   \label{inner_porfolio_constraint1}  \\
&\Pi_k\geq0, ~~C_k\geq 0, \label{inner_porfolio_constraint2} \\
&C_k \leq R_f(W_k- 1^{\top}_n\Pi_k) , ~ \text{for}~ k=0,\cdots,K-1. \label{inner_porfolio_constraint3}
\end{align}
\end{subequations}
Note that (\ref{inner_porfolio_constraint4})-(\ref{inner_porfolio_constraint1}) are equivalent to (\ref{discrt_market_state})-(\ref{discrt_wealth}), and(\ref{inner_porfolio_constraint2})-(\ref{inner_porfolio_constraint3}) are equivalent to (\ref{discrt_constraint}). The advantage of this reformulation is that the inner optimization problem (\ref{inner_portfolio_opt}) has linear constraints. Therefore, we may find the global maximizer of (\ref{inner_portfolio_opt}) as long as the objective function in (\ref{inner_porfolio_obj}) is jointly concave in $\textbf{a}$.

Heuristically,  we need to design near-optimal penalty functions in order to  obtain tight dual bounds on $H_0$. A natural approach is to investigate the optimal penalty $M^*$ for the discrete-time problem according to (\ref{opt_penalty_MDPs}):
$$M^{*}(\textbf{a},\textbf{Z},\tilde{\textbf{Z}})=\sum_{k=0}^{K-1}\Delta H_{k+1}(\textbf{a},\textbf{Z},\tilde{\textbf{Z}}),$$
where $\Delta H_{k+1}$ is the deviation in $H_{k+1}$ from the conditional mean. In practice we can approximate $H_{k}$ by $\hat{H}_{k}=W_k^{1-\gamma}\hat{J}_{k}$; however, it does not mean that $\Delta \hat{H}_{k+1}$ can be easily computed, since an intractable conditional expectation (that is, $\mathbb{E}_k[\hat{H}_{k+1}]$) over $(n+d)$-dimensional space  is involved.  Another difficulty is that  $M^{*}=\sum_{k=0}^{K-1}\Delta H_{k+1}$  enters  into (\ref{inner_porfolio_obj}) with possibly positive or negative signs for different realizations of $(\textbf{Z},\tilde{\textbf{Z}})$, making  the objective function of (\ref{inner_portfolio_opt}) nonconcave,  even if $U$ is a concave function. Therefore, it might be extremely hard to locate the global maximizer of (\ref{inner_portfolio_opt}).

To address these problems, we exploit the value function-based optimal penalty $h^{*}_v$  for the continuous-time problem  (\ref{continuous_market_state})-(\ref{Value_portfolio}), recalling that our discrete-time problem is discretized from the continuous-time model. Based on the form of   $h^{*}_v$  we will propose \emph{a dual feasible penalty in the sense of  (\ref{portfolio_penalty_set}) for the discrete-time problem,} which is also easy to compute. Assuming that all the technical conditions in Theorem\ref{Theorem5:StrongDual_IdealPenalty} hold, we  can apply the result (\ref{optimal_penalty_2}) by selecting $x_t=(\phi_t,W_t)$, $V(t,x_t)=
V(t,\phi_t,W_t)$, $\sigma(t,x_t)=\begin{pmatrix}
\sigma^{\phi,1}_t &\sigma^{\phi,2}_t\\
W_t\pi_t\sigma_t &0
\end{pmatrix}$, and  $dw_t=\begin{pmatrix}
dz_t   \notag\\
d\tilde{z}_t
\end{pmatrix}$ such that
\begin{align}
h_v^*(\textbf{u},\textbf{z},\tilde{\textbf{z}})=&\int_{0}^{T}
\begin{pmatrix}
V_{\phi}(t,\phi_t,W_t) \\
V_{W}(t,\phi_t,W_t)
\end{pmatrix}^{\top}
\begin{pmatrix}
\sigma^{\phi,1}_t &\sigma^{\phi,2}_t\\
W_t\pi_t\sigma_t &0
\end{pmatrix}
\begin{pmatrix}
dz_t   \notag\\
d\tilde{z}_t
\end{pmatrix} \notag\\
=&\sum_{k=0}^{K-1}\int_{k\delta}^{(k+1)\delta}\bigg[V^{\top}_{\phi}(t,\phi_t,W_t)\sigma^{\phi,1}_t dz_t  \notag\\
&+V^{\top}_{\phi}(t,\phi_t,W_t)\sigma^{\phi,2}_t d\tilde{z}_t  +V_{W}(t,\phi_t,W_t)W_t\pi_t\sigma_tdz_t\bigg]  \notag\\
=&\sum_{k=0}^{K-1}\int_{k\delta}^{(k+1)\delta}\beta^t \bigg[ W_t^{1-\gamma}\nabla_{\phi}\tilde{J}^{\top}(t,\phi_t)\sigma^{\phi,1}_tdz_t + W_t^{1-\gamma}\nabla_{\phi}\tilde{J}^{\top}(t,\phi_t)\sigma^{\phi,2}_td\tilde{z}_{t}   \notag\\
&+(1-\gamma)W_t^{1-\gamma}\tilde{J}(t,\phi_t)\pi_t\sigma_{t}dz_t\bigg],  \label{continuous_portfolio_martingale}
\end{align}
for  $\textbf{u}=(\pi_t,\tilde{c}_t)_{0\leq t \leq T}\in \mathcal{U}_{\mathbb{F}}(0)$, and the last equality holds due to the structure of the value function (\ref{value_function_SDE_structure}). In particular, we use $\nabla_{\phi}\tilde{J}$ to denote the gradient of the function $\tilde{J}$ with respect to $\phi$. By discretizing the Ito stochastic integrals  in (\ref{continuous_portfolio_martingale}), we propose a heuristic -- using the $(k+1)$-th term in the summation -- to approximate  $\Delta H_{k+1}$ in $M^*$, that is,
\begin{align}
 \Delta H_{k+1}\approx &\beta^{ k\delta}
\big[ W_k^{1-\gamma}\nabla_{\phi}J^{\top}_{k}(\phi_k)\sigma^{\phi,1}_k\sqrt{\delta}Z_{k+1} \notag\\
&+ W_k^{1-\gamma}\nabla_{\phi}J^{\top}_{k}(\phi_k)\sigma^{\phi,2}_k\sqrt{\delta}\tilde{Z}_{k+1} \notag\\
&+(1-\gamma)W_k^{-\gamma}J_{k}(\phi_k)\Pi_k^{\top}\sigma_k\sqrt{\delta}Z_{k+1}\big], \label{penalty_portfolio0}
\end{align}
where we use $J_{k}(\phi)$ to approximate $\tilde{J}(k\delta,\phi)$ and also use the substitution $\Pi_k=W_k\pi_k$.

We then describe a procedure to  empirically approximate $M^*=\sum_{k=0}^{K-1}\Delta H_{k+1}$ based on (\ref{penalty_portfolio0})  using simulation. Given a realization of $(\textbf{Z},\tilde{\textbf{Z}})$  we can obtain the realized terms of $\bar{\phi}_k\triangleq \phi_k(\phi_0,\textbf{Z},\tilde{\textbf{Z}})$, $\bar{\sigma}_k\triangleq\sigma(\bar{\phi}_k)$, $\bar{\sigma}^{\phi,1}_{k}\triangleq\sigma^{\phi,1}(k,\bar{\phi}_k)$, $\bar{\sigma}^{\phi,2}_k\triangleq\sigma^{\phi,2}(k,\bar{\phi}_k)$; with an admissible strategy $\hat{\textbf{a}}=(\hat{a}_0,\cdots,\hat{a}_{K})$, we can also  obtain
 $\bar{W}_k\triangleq W_k(W_0,\hat{\textbf{a}}(\phi_0,W_0,\textbf{Z}_k,\tilde{\textbf{Z}}_k),\textbf{Z}_k,\tilde{\textbf{Z}}_k)$ via (\ref{discrt_wealth}) as an approximation to the wealth under the optimal policy. Then we can approximate $M^{*}(\textbf{a},\textbf{Z},\tilde{\textbf{Z}})$  by
 \begin{align}
M_1(\textbf{a},\textbf{Z},\tilde{\textbf{Z}})\triangleq\sum_{k=0}^{K-1}\left(\Psi^{1}_{k}\left(\textbf{a},\textbf{Z},\tilde{\textbf{Z}}\right) Z_{k+1}+\Psi^{2}_{k}(\textbf{a},\textbf{Z},\tilde{\textbf{Z}}) \tilde{Z}_{k+1}\right), \label{MarkovDiff:penalty_portfolio1}
\end{align}
where
\begin{align}
\Psi^{1}_{k}(\textbf{a},\textbf{Z},\tilde{\textbf{Z}})=&\beta^{ k\delta}
\left[\bar{W}_k^{1-\gamma}\Xi^{2\top}_{k}\left(\bar{\phi}_k\right)\bar{\sigma}^{\phi,1}_k\sqrt{\delta} + (1-\gamma)\bar{W}_k^{-\gamma}\Xi^1_{k}(\bar{\phi}_k)\Pi_k^{\top}\bar{\sigma}_k\sqrt{\delta}\right],   \label{MarkovDiff:integrand2}\\
\Psi^{2}_{k}(\textbf{a},\textbf{Z},\tilde{\textbf{Z}})=&\beta^{ k\delta}\bar{W}_k^{1-\gamma}\Xi^{2\top}_{k}(\bar{\phi}_k)\bar{\sigma}^{\phi,2}_k\sqrt{\delta}, \notag
\end{align} and where   $\Xi^1_{k}(\cdot)$  is a scalar function of $\phi$, whereas $\Xi^2_{k}(\cdot)$ is an $m$-dimensional  function  of $\phi$. As suggested by  (\ref{penalty_portfolio0}), $\Xi^1_{k}(\cdot)$ and $\Xi^2_{k}(\cdot)$  are preferably chosen as  $\hat{J}_k(\cdot)$ -- an approximation of $J_k(\cdot)$, and $\nabla_{\phi}\hat{J}_{k}(\phi_k)$ -- an approximation of  $\nabla_{\phi}J_{k}(\phi_k)$, respectively. In the case that $\hat{J}_{k}(\phi)$ is not differentiable in $\phi$, we may apply the the finite difference method on $\hat{J}_{k}(\phi_k)$ to obtain the difference quotient as $\Xi^2_{k}(\cdot)$ (i.e., a nominal approximation of $\nabla_{\phi}\hat{J}_{k}(\phi_k)$). We  verify  in Proposition \ref{Theorem6:Valid_UP_BOUND1} below that $M_1$ is dual feasible and hence $\mathcal{L}M_1$ is an upper bound on $H_0$.

It remains to show why the forms of $\Psi^{1}_{k}$ and $\Psi^{2}_{k}$ make the inner optimization problem (\ref{inner_portfolio_opt}) easy to solve. This is because both functions are \emph{affine} in $\textbf{a}$, regardless of the realizations of $\textbf{Z}$ and $\tilde{\textbf{Z}}$. To be specific, when a realization of $(\textbf{Z},\tilde{\textbf{Z}})$ is fixed, $\Psi^{2}_{k}$ is a constant with respect to $\textbf{a}$, while $\Psi^{1}_{k}$ is affine in  $\Pi_k$ (hence, in \textbf{a}).  Therefore, together with the concave property of $U(\cdot)$, the inner optimization problem (\ref{inner_portfolio_opt}) is guaranteed to be convex   with $M=M_1$. To find some variants of the penalties while still keeping the convexity of the inner optimization problem, we also generate $\breve{\Psi}^{1}_{k+1}$ based on a first-order Taylor expansion of  $\Psi^{1}_{k+1}$ in (\ref{MarkovDiff:integrand2}) around the strategy $\hat{a}_{k-1}$, $k=1,\cdots,K$ (we only expand the first term, since the second term is already linear in $\Pi_k$):
\begin{align*}
\breve{\Psi}^{1}_{k+1}(\textbf{a},\textbf{Z},\tilde{\textbf{Z}})= &\beta^{ k\delta}\big[\bar{W}_k^{1-\gamma}+ (1-\gamma) \bar{W}_k^{-\gamma} \notag\big((\bar{R}_{k}-R_f\textbf{1}_n)^{\top}(\Pi_{k-1}-\bar{\Pi}_{k-1}) \notag\\
&-(C_{k-1}-\bar{C}_{k-1})\big)\big]\cdot\Xi^{2\top}_{k}(\bar{\phi}_k)\bar{\sigma}^{\phi,1}_k\sqrt{\delta}+ \beta^{ k\delta}(1-\gamma)\bar{W}_k^{1-\gamma}\Xi^1_{k}(\bar{\phi}_k)\Pi_k^{\top}\bar{\sigma}_k\sqrt{\delta}, %\label{penalty_portfolio2}
\end{align*}
where $\bar{R}_{k}\triangleq R_{k}(\phi_0,\textbf{Z},\tilde{\textbf{Z}})$, $(\bar{\Pi}_{k},\bar{C}_{k})\triangleq \hat{a}_{k}(\phi_0,W_0,\textbf{Z},\tilde{\textbf{Z}})$. Then $\breve{\Psi}^{1}_{k+1}$ is affine in $\Pi_{k-1}$ and $C_{k-1}$. We can also obtain a variant of $\Psi^{2}_{k+1}$ that is is affine in $\Pi_{k-1}$ and $C_{k-1}$,  say $\breve{\Psi}^{2}_{k+1}$, in exactly the same way. In our numerical examples we will  consider dual bounds generated by $M_1$ as well as $M_2$, where
\begin{align}
M_2(\textbf{a},\textbf{Z},\tilde{\textbf{Z}})\triangleq\sum_{k=0}^{K-1}\left(\breve{\Psi}^{1}_{k}(\textbf{a},\textbf{Z},\tilde{\textbf{Z}}) Z_{k+1}+\breve{\Psi}^{2}_{k}(\textbf{a},\textbf{Z},\tilde{\textbf{Z}}) \tilde{Z}_{k+1}\right). \label{MarkovDiff:penalty_portfolio3}
\end{align}

To go further, we can also generate a penalty function by linearizing $\Psi^{1}_{k+1}$ around $(\hat{a}_0,\cdots,\hat{a}_{k-1})$.  We show  $M_2\in \mathcal{M}_\mathbb{G}(0)$ in Proposition \ref{Theorem6:Valid_UP_BOUND1} as well.

\begin{Proposition} \label{Theorem6:Valid_UP_BOUND1}
Both $M_1$ and $M_2$ are dual feasible in the sense of  (\ref{portfolio_penalty_set}),  i.e., $M_1,M_2\in \mathcal{M}_\mathbb{G}(0)$. Hence, both $\mathcal{L}M_1$ and $\mathcal{L}M_2$ are upper bounds on $H_{0}$.
\end{Proposition}
\begin{IEEEproof}
 First, we  show that $\Psi^{i}_{k}(\textbf{a},\textbf{Z},\tilde{\textbf{Z}})$ is   $\mathcal{G}_k$-adapted  given any $\textbf{a}\in \mathbb{A}_{\mathbb{G}}$ for $i=1,2$. Noting   that $\bar{\phi}_k$, $\Xi^{1}_{k}(\bar{\phi}_k)$, $\Xi^{2}_{k}(\bar{\phi}_k)$, $\bar{\sigma}_{k}$, $\bar{\sigma}^{\phi,j}_{k}$ ($j=1,2$), and $\bar{W}_k$ are naturally $\mathcal{G}_k$-adapted under a fixed non-anticipative policy $\hat{\textbf{a}}\in \mathbb{A}_{\mathbb{G}}$. Therefore, $\Psi^{2}_{k+1}(\textbf{a},\textbf{Z},\tilde{\textbf{Z}})$ is  $\mathcal{G}_k$-adapted.  We also observe that  $\Pi_k$  is  $\mathcal{G}_k$-adapted as $\textbf{a}\in \mathbb{A}_{\mathbb{G}}$; therefore,  $\Psi^{1}_{k}(\textbf{a},\textbf{Z},\tilde{\textbf{Z}})$ is  $\mathcal{G}_k$-adapted for any  $\textbf{a}\in \mathbb{A}_{\mathbb{G}}$.

Second, since $Z_{k+1}$ and $\tilde{Z}_{k+1}$ have zero means and  are independent of $\mathcal{G}_k$ and $(\phi_0,W_0)$,  along with the linearity of  $\Psi^{1}_{k}$ (resp., $\Psi^{2}_{k}$) in  $Z_{k+1}$ (resp., $\tilde{Z}_{k+1})$, we have for  $k=0,\cdots,K-1$,
\begin{align*}
\mathbb{E}\left[\Psi^{1}_{k}\cdot Z_{k+1}\big|\phi_0,W_0\right]&=\mathbb{E}_0\left[\Psi^{1}_{k}\cdot \mathbb{E}_k[Z_{k+1}\big]\right]=0 \text{~~for all~~} \textbf{a}\in \mathbb{A}_{\mathbb{G}};\\
\mathbb{E}\left[\Psi^{2}_{k}\cdot\tilde{Z}_{k+1}\big|\phi_0,W_0\right]&=\mathbb{E}_0\left[\Psi^{2}_{k}\cdot\mathbb{E}_k[\tilde{Z}_{k+1}]\right]=0 \text{~~for all~~} \textbf{a}\in \mathbb{A}_{\mathbb{G}}.
\end{align*}
Therefore, $ \mathbb{E}[M_1(\textbf{a},\textbf{Z},\tilde{\textbf{Z}})|\phi_0,W_0]= 0$ for all $\textbf{a}\in \mathbb{A}_{\mathbb{G}},$ and hence $M\in \mathcal{M}_\mathbb{G}(0)$. The same argument can also apply on $M_2$. Therefore, $M_2\in \mathcal{M}_\mathbb{G}(0)$.
\end{IEEEproof}

The penalties in the form of (\ref{MarkovDiff:penalty_portfolio1}) or (\ref{MarkovDiff:penalty_portfolio3}) bear several advantages. First, unlike $\sum_{k=0}^{K-1}\Delta \hat{H}_{k+1}$ that directly approximates the optimal penalty of the discrete-time model,  our proposed penalties (\ref{MarkovDiff:penalty_portfolio1}) and (\ref{MarkovDiff:penalty_portfolio3}) does not involve  any conditional expectation and can be evaluated efficiently; therefore, a substantial amount of computational work can be avoided. Second, the design of such penalties is quite flexible: we can use any admissible policy to obtain a dual feasible penalty, and linearize around this policy if necessary, which guarantees the convexity of  the inner optimization problem (\ref{inner_portfolio_opt}).

%\begin{Remark}
%For the continuous-time model we should verify that $V(t,\phi_t,W_t)$ satisfies certain conditions such that (\ref{continous_portfolio_martingale}) is a martingale and hence it is a valid penalty. In the discrete-time model we do not have this verification step, since we can decompose each term in (\ref{penalty_portfolio1}) as a product of $Z_{k+1}(\tilde{Z}_{k+1})$ and a $\mathcal{G}_k$-adapted term so that it is automatically a martingale term; further, the differentiability of $J_k(\phi_k)$ does not influence the validity of (\ref{penalty_portfolio1}) and (\ref{penalty_portfolio3}) being  penalty functions, either.
%\end{Remark}

%We should point out that  the application of the dual approach on continuous-state MDPs is mainly limited by the tractability of the inner optimization problem, especially when its feasible region  that is constrained by the state equation and the control set (for example, (\ref{inner_porfolio_constraint1})-(\ref{inner_porfolio_constraint3})) is not convex.  However, with our proposed penalty the inner optimization problem is always tractable, as long as the inner optimization problem with zero penalty $(i.e., M(\textbf{a},\textbf{Z},\tilde{\textbf{Z}})=0)$ can be solved.

\subsection{Numerical Experiments}
In this section we discuss the use of Monte Carlo simulation to evaluate the performance of the suboptimal policies and the dual bounds on the expected utility (\ref{discrt_value}).  We consider a model with three risky assets $(n=3)$ and one market state variable $(m=1)$.  The dynamics (\ref{continuous_asset})-(\ref{continuous_market_state}) of the market state and assets returns are the same as those  considered in \cite{haugh:2006}. In particular, let $\mu^{\phi}_k=-\lambda\phi_{k}$, $\mu_{k}=\mu_0+\mu_1 \phi_k $, $\sigma_k\equiv \sigma$, $\sigma^{\phi,1}_k\equiv \sigma^{\phi,1}$, and $\sigma^{\phi,2}_k\equiv  \sigma^{\phi,2}$, in (\ref{discrt_market_state})-(\ref{discrt_R}).  The parameter values are listed in the following tables including $r_f$, $\lambda$, $\mu_0$, $\mu_1$, $\sigma$, $\sigma^{\phi,1}$, and $\sigma^{\phi,2}$. Note from  (\ref{continuous_market_state}) that the market state $\phi$ follows a mean-reverting Ornstein-Uhlenbeck process: it has relatively small mean reversion rate and volatility in the parameter set 1, while it has relatively large mean reversion rate and volatility in the parameter set 2.  We choose $T=1$ year and $\delta=0.1$ year in our numerical experiments. In addition, we use  $\alpha=0.5$  for the weight of the intermediate utility function and use $\beta=1$ as the discount factor. We assume $\phi_0=0$ and $W_0=1$ as the initial condition and impose the constraint (\ref{discrt_constraint}) on the control space $\mathcal{A}$ in the following numerical tests.
%\begin{sideways}
%\begin{landscape}
\begin{table}[htpb]
  \centering \caption{Parameter Set 1 }\label{Parameter1}
  \setlength{\tabcolsep}{0.8pt}

\begin{tabular}{|c|cc|c|c|} \hline

    %$\pi$  &  $\theta$ \\
    $$ &$\mu_0$ &$\mu_1$  &$\sigma$  &$r_f$\\
   \cline{2-5}
        %\rule[1ex]{0pt}{1em}\\

  $\log (R)$
  &$\left(
    \begin{array}{cc}
        $0.081$ \\
        $0.110$ \\
        $0.130$ \\
        \end{array}
    \right)$

   & $\left(
        \begin{array}{c}
        $0.034$\\
        $0.059$\\
        $0.073$\\
        \end{array}
          \right)$

  & $\left(
    \begin{array}{ccc}
     $ 0.186$  &$0.000$   &$0.000$ \\
     $0.228$  &$0.083$     &$0.000$ \\
    $0.251$     &$0.139$       &$0.069$ \\
    \end{array}
  \right)$

  & 0.01\\
  \hline
  \hline
  $\phi$ &\multicolumn{2}{c|} {$\lambda$ } &$\sigma^{\phi,1}$ &$\sigma^{\phi,2}$\\
  \cline{2-5}
  $$
  &\multicolumn{2}{c|} {0.336 }
  & $\left(
    \begin{array}{ccc}
     $-0.741$     &$-0.037$       &$-0.060$  \\
  \end{array}
  \right)$
  &$0.284$ \\
  \hline
\end{tabular}
  \end{table}

\begin{table}[htpb]
  \centering \caption{Parameter Set 2  }\label{Parameter2}
  \setlength{\tabcolsep}{0.8pt}
  \begin{tabular}{|c|cc|c|c|} \hline

    %$\pi$  &  $\theta$ \\
    $$ &$\mu_0$ &$\mu_1$  &$\sigma$  &$r_f$\\
   \cline{2-5}
        %\rule[1ex]{0pt}{1em}\\

  $\log (R)$
  &$\left(
    \begin{array}{cc}
        $0.081$ \\
        $0.110$ \\
        $0.130$ \\
        \end{array}
    \right)$

   & $\left(
        \begin{array}{c}
        $0.034$\\
        $0.059$\\
        $0.073$\\
        \end{array}
          \right)$

  & $\left(
    \begin{array}{ccc}
     $ 0.186$  &$0.000$   &$0.000$ \\
     $0.228$  &$0.083$     &$0.000$ \\
    $0.251$     &$0.139$       &$0.069$ \\
    \end{array}
  \right)$

  & 0.01\\
  \hline
  \hline
  $\phi$ &\multicolumn{2}{c|} {$\lambda$ } &$\sigma^{\phi,1}$ &$\sigma^{\phi,2}$\\
  \cline{2-5}
  $$
  &\multicolumn{2}{c|} {1.671}
  & $\left(
    \begin{array}{ccc}
     $-0.017$     &$0.149$       &$-0.058$  \\
  \end{array}
  \right)$
  &$1.725$ \\
  \hline
\end{tabular}
  \end{table}

For each parameter set we first use the discrete state-space approximation method to solve the recursion (\ref{recursion_portfolio}). In particular, we approximate the market state variable  $\phi_k$ using a grid with $21$ equally spaced grids from $-2$ to $2$, and the transition between these grid points is determined by (\ref{discrt_market_state}) noting that $\phi_{k+1}\sim N\big(\phi_k+\mu^{\phi}_k \delta, (\parallel \sigma^{\phi,1}_k \parallel^2+ \parallel \sigma^{\phi,2}_k \parallel^2)\delta \big) $; the random variables $Z_{k}$ and $\tilde{Z}_{k}$ are approximated by  Gaussian quadrature method with $3$ points for each dimension (see, e.g., \cite{judd:1998}). So the joint distribution of the market state and the returns  are approximated by  a total of $3^3\times 21=567$ grid points, which are used to compute the conditional expectation in  (\ref{recursion_portfolio}): we assume $\phi_{k+1}$ and $R_{k+1}$  are independent conditioned on $\phi_k$, then the conditional expectation reduces to a finite weighted sum. For the optimization problem in (\ref{recursion_portfolio}) we use CVX (\cite{Grant2013}), a package to solve convex optimization problems in MATLAB, to determine the optimal consumption and investment policy on each grid of $\phi_k$ at time $k$. We record the value function and the corresponding policy on this grid at each time $k=0,\cdots, K.$ Note that the market state variable $\phi_k$ is one dimensional, so the value function and the policy can be naturally defined on the market state $\phi_k$ that is outside the grid  by piecewise linear interpolation.  In our numerical implementation  the extended value function and the extended policy play the roles of  $\Xi^{1}_{k}(\phi)$ (i.e., $\hat{J}_k(\phi)$) and the approximate policy $\hat{\textbf{a}}$ to the discrete-time  problem (\ref{discrt_model})-(\ref{discrt_value}); and we  take the slope of the piecewise linear function $\Xi^{1}_{k}(\phi)$ as $\Xi^{2}_{k}(\phi)$, if $\phi$ is between the grid points; otherwise, we can use the average slope of two consecutive lines as $\Xi^{2}_{k}(\phi)$.
% and the transition  between these grid points  is determined to match the mean and covariance information contained in (\ref{discrt_market_state})-(\ref{discrt_R}).

We then repeatedly  generate random sequences of $(\textbf{Z},\tilde{\textbf{Z}})$, based on which we generate  the sequences of market states and returns  according to their joint probability distribution  (\ref{discrt_model})-(\ref{discrt_value}). Then we apply the aforementioned policy  $\hat{\textbf{a}}$  on these sequences to get an estimate of the lower bound on the value function $H_{0}$; based on each random sequence we can also solve the inner optimization problem (\ref{inner_portfolio_opt}) with  penalty $M_1$  in (\ref{MarkovDiff:penalty_portfolio1}) or $M_2$ in  (\ref{MarkovDiff:penalty_portfolio3}), which leads to an estimate of the upper bound on $H_0$. We present our numerical results in the following tables:  the lower bound, which is referred to as ``Lower Bound'', is obtained by generating $100$ random sequences of $(\textbf{Z},\tilde{\textbf{Z}})$ and their antithetic pairs (see \cite{glasserman:2004} for an introduction on antithetic variates) in a single run and a total number of $10$ runs; the upper bounds induced by penalties $M_1$ and $M_2$, which are referred to as ``Dual Bound 1'' and ``Dual Bound 2'' respectively, are obtained by generating $30$ random sequences of $(\textbf{Z},\tilde{\textbf{Z}})$ and their antithetic pairs  in a single run and a total number of $10$ runs. To see the effectiveness of these proposed penalties,  we use zero penalty and  repeat the same procedure to  compute the upper bounds that are referred to as ``Zero Penalty'' in the table. These bounds on the value function $H_0$ (i.e., the expected utility) are reported in the sub-column ``Value'', where each entry shows the sample average and the standard error (in parentheses) of the $10$ independent runs. We also compute the certainty equivalent of the  expected utility in the sub-column ``CE'' (this is reported in the literature such as \cite{cvitanic:2003}), where ``CE'' is defined through $U(\text{CE})=\text{Value}.$ For ease of comparison, we compute the duality gaps --  differences of the lower bound with each upper bound on the expected utility and its certainty equivalent -- as a fraction of the lower bounds, and list the smaller fraction in the column ``Duality Gap''.

{\footnotesize

\begin{table*}[t]
\centering \caption{ Results with Parameter Set $1$ } \label{Table:Result1}
 \begin{tabular}{|c|cc|cc|cc|cc|cc|} \hline
  \multicolumn{1}{|c} {} &\multicolumn{2}{|c|} {Lower Bound}  &\multicolumn{2}{c|} {Dual Bound 1} &\multicolumn{2}{c|} { Dual Bound 2}   &\multicolumn{2}{c|} { Zero Penalty } &\multicolumn{2}{c|} { Duality Gap }\\ \hline
    { $\gamma$ }  &Value &CE &Value &CE   &Value  &CE &Value  &CE   &Value  &CE
   \rule[1ex]{0pt}{1em}\\\hline\hline
    $1.5$       &$-5.480$    &$0.1332$      &$ -5.391$   &$0.1376$         &$-5.392$     &$0.1376 $    &-4.861   & 0.1693   &1.61\% &3.30\% \\
    $$          &$(0.003)$   &$(0.0001)$    &$(0.008)$   &$(0.0004)$       &$(0.007)$    &$(0.0004)$   &(0.012)  &(0.0008)  &       & \\
    $3.0$       &$-42.887$   &$0.1080$     &$-39.227$    &$0.1129 $        &$-39.873$    &$0.1120 $    &-27.562  &0.1347    &7.53\% & 3.70\% \\
    $$          &$(0.036)$   &$(0.0001)$    &$(0.164)$   &$(0.0002)$       &$(0.317)$    &$(0.0004)$   &(0.252)  &(0.0006)  &       & \\
    $5.0$       &$-2445.9$   &$0.1005$     &$-2066.5$    &$0.1049$         &$-2025.5$    &$0.1054$     &-1105.7  &0.1226    &15.51\% &4.38\% \\
    $$          &$(1.635)$   &$(0.0001)$    &$(22.019)$  &$(0.0003)$       &$(17.833)$   &$(0.0002)$   &(16.438) &(0.0004)  & &  \\

    \hline
\end{tabular}
\end{table*}

\begin{table*}[t]
\centering \caption{ Results with Parameter Set $2$} \label{Table:Result2}
  \begin{tabular}{|c|cc|cc|cc|cc|cc|} \hline
  \multicolumn{1}{|c} {} &\multicolumn{2}{|c|} {Lower Bound} &\multicolumn{2}{c|} {Dual Bound 1} &\multicolumn{2}{c|} {Dual Bound 2} &\multicolumn{2}{c|} { Zero Penalty } &\multicolumn{2}{c|} { Duality Gap }\\ \hline
    { $\gamma$ }  &Value &CE &Value &CE    &Value  &CE &Value  &CE  &Value  &CE
   %{ $Y_{0}$    } &$\text{Lower Bound}$ &$\text{\emph{l}}$  &$m=100$ &$m=500$ &$m=1000$

    \rule[1ex]{0pt}{1em}\\\hline\hline
    $1.5$      &$-5.466$      &$0.1339$    &$ -5.380$   &$0.1382 $     &$-5.381$     &$ 0.1381$    &-4.864   &0.1691    &1.56\% &3.14\%  \\
    $$         &$(0.005)$     &$(0.0001)$  &$(0.011)$   &$(0.0006)$    &$(0.015)$    &$(0.0008)$   &(0.020)  &(0.0008)  & & \\
    $3.0$      &$-42.585$     &0.1084      &$-39.645$   &$0.1123$      &$-39.690$    &$0.1122$     &-27.708  &0.1343    &6.80\% &3.51\%  \\
    $$         &$(0.081)$     &(0.0001)    &$(0.229)$   &$(0.0003)$    &$(0.155)$    &$(0.0002)$   &(0.209)  &(0.0005)  & &\\
    $5.0$      &$-2431.6$     &$0.1007$    &$-2043.8$   &$ 0.1052 $    &$-2040.7$    &$0.1052 $    &-1122.1  &0.1222    &15.95\% &4.47\%  \\
    $$         &$(7.510)$     &(0.0001)    &$(11.881)$  &$(0.0002)$    &$(19.882)$   &$(0.0003)$   &(9.842)  &(0.0004)  & & \\

    \hline
\end{tabular}
\end{table*}}

We consider utility functions with different relative risk aversion coefficients $\gamma=1.5,3.0,$ and $5.0$, which reflect low, medium and high degrees of risk aversions.
The  dual bounds induced by zero penalty perform poorly as we expected. On the other hand, it is hard to distinguish the performance of ``Dual Bound 1'' and ``Dual Bound 2'', which may imply that the second term in   (\ref{MarkovDiff:integrand2}) plays an essential role in the inner optimization problem in order to make the dual bounds tight in this problem. We observe that the duality gaps on the value function $H_0$ are generally smaller when $\gamma$ is small, implying that both the approximate policy and penalties are near optimal. For example, when $\gamma=1.5$, the duality gaps are within $2\%$  of the optimal expected utility  for all sets of parameters. As $\gamma$ increases, the duality gaps generally become larger.

There are several possible reasons for the enlarged duality gaps on the value function  with increasing $\gamma$. Note that the utility function $U(x)$ is a power function (with negative power of $1-\gamma$) of $x$  and it decreases at a higher rate with larger $\gamma$, as $x$ approaches zero.  This is  reflected by the fact that both the lower and upper bounds on the value function $H_0$ decrease rapidly  with higher value of $\gamma$. In the case of evaluating the upper bounds on $H_0$, it can be inferred that with larger $\gamma$ the objective value (\ref{inner_porfolio_obj})  is more sensitive to the solution of the inner optimization problem (\ref{inner_portfolio_opt}), and hence the quality of the penalty functions. In other words, even a small torsion of the optimal penalty will lead to a significant deviation of the dual bound. In our case  the heuristic penalty is derived by discretizing the value function-based penalty for the continuous-time problem; however, this penalty may become  far away from optimal for the discrete-time problem when $\gamma$ increases. Similarly, obtaining tight lower bounds on the expected utility by simulation under a sub-optimal policy also suffers the same problem, that is, solving a sub-optimal policy based on the same approximation scheme of the recursion (\ref{recursion_portfolio}) may cause more utility loss  with larger $\gamma$. The performance of the sub-optimal policy also influences the quality of the penalty function, since the penalties $M_1$ and $M_2$ involve the wealth $\bar{W}_k$ induced by the suboptimal policy and its error compared with the wealth under the optimal policy will be accumulated over time. Hence, the increasing duality gaps on the value function with larger risk aversion coefficients  are contributed by  both sub-optimal policies and sub-optimal penalties.

 These numerical results provide us with some guidance in terms of computation when we apply  the dual approach: we should be more careful with  designing the penalty function if the objective value of the inner optimization problem is numerically sensitive either to its optimal solution or to the choice of the penalty function.   Fortunately, the sensitivity of the expected utility  with respect to $\gamma$ in this problem is relieved to some extent by considering its certainty equivalent. We can see from the table that the differences between the lower bounds and the upper bounds in terms of  ``CE'' are kept at a relatively constant range for different values of $\gamma$.

%To improve the quality of the the dual bound,  we need  more computational work to capture a better suboptimal policy, which implicitly implies

\section{Conclusion}\label{Section:Conclusion}
In this paper we study the dual formulation of controlled Markov diffusions by means of information relaxation. This dual formulation provides new insights into seeking the value function: if we can find an optimal solution to the dual problem, i.e., an optimal penalty, then the value function can be recovered without solving the HJB equation. From a more practical point of view, this dual formulation can be used to find a  dual bound on the value function. We explore the structure of the value function-based optimal penalty, which provides the theoretical basis for developing near-optimal penalties that lead to tight dual bounds. As in the case of MDPs, if we compare the dual bound on the value function of a controlled Markov diffusion with the lower bound generated  by Monte Carlo simulation under a sub-optimal policy, the duality gap can serve as an  indication on how well the sub-optimal policy performs and how much we can improve on our current policy. Furthermore, we also expose the connection of the  gradient-based optimal penalty between controlled Markov diffusions and MDPs in Appendix.

We carried out numerical studies in a dynamic portfolio choice problem that is discretized from a continuous-time model. To derive tight dual bounds on the expected utility, we proposed a class of penalties that can be viewed as discretizing the value function-based optimal penalty of the continuous-time problem, and these new penalties make the inner optimization problem  computationally tractable. This approach  has potential use in many other interesting applications where the system dynamic is modeled as a controlled Markov diffusion. Moreover, we investigate the sensitivity of the quality of both lower and upper bounds in terms of duality gaps with respect to different parameters. These numerical studies complement the existing examples of applying the dual approach to continuous-state MDPs.

\bibliographystyle{IEEEtran}
\bibliography{Zhou-Bibtex}

%\begin{thebibliography}{99}

%\end{thebibliography}

\appendix%{Complement of Section \ref{Section:Dual_Diffusion}}\label{Section:Appendix:Section2.1}

%\section{Complement of Section \ref{Section:Dual_Diffusion}}\label{Section:Appendix:Section2.1}
In the appendix we aim to develop the value function-based penalty as a solution to the dual problem on the right side of (\ref{strong_duality_1}), which can be viewed as  the counterpart of (\ref{opt_penalty_MDPs}) in the setting of controlled Markov diffusions. For this purpose we need to define a solution to  the stochastic differential equation(SDE) (\ref{state_SDE}) with an anticipative control $\textbf{u}\in \mathcal{U}(0)$. Therefore, we introduce the Stratonovich calculus and anticipating stochastic differential equation in Appendix \ref{Section:Appendix:anticipating}, and present the value function-based optimal penalty in Appendix \ref{Section:Appendix:value_penalty}.   We also review the dual representation of the optimal stopping problem under the diffusion process in Appendix \ref{Section:Appendix:optimal_stopping}.

\vspace{-2mm}
\subsection{Anticipating Stochastic Differential Equation}\label{Section:Appendix:anticipating}

%The SDE is in Ito sense with an admissable strategy u.
 There are several ways to integrate stochastic processes that are not adapted to Brownian motions such as Skorohod and (generalized) Stratonovich integrals (see, e.g, \cite{nualart2006malliavin,Ocone1989}).  In this subsection we present the Stratonovich integral and its associated Ito formula.  Then we generalize  the controlled diffusion (\ref{state_SDE}) to the Stratonovich sense following \cite{davis1992deterministic}.

We first assume that $\textbf{w}=(w_t)_{t \in[0,T]}$ is a one-dimensional Brownian Motion in the probability space $(\Omega,\mathcal{ F},\mathbb{P})$. We denote by $I$ an arbitrary partition of the interval $[0,T]$ of the form $I=\{0= t_0 < t_1<\cdots <t_n=T\}$% and define $|I|=\sup_{0\leq i\leq n-1}(t_{i+1}-t_i)$.
\begin{definition}(Definition 3.1.1 in \cite{nualart2006malliavin})\label{Stratonovich_integral}
 We say that a measurable process $\textbf{y}=(y_t)_{t \in[0,T]}$ such that $\int_0^T|y_t|dt<\infty$ a.s. is Stratonovich integrable if the family
 $$S^{I}=\int_0^T y_t \sum_{i=0}^{n-1}\frac{w_{t_i+1}-w_{t_i}}{t_{i+1}-t_i}\mathbbm{1}_{(t_i,t_{i+1}]}(t)\, dt$$
 converges in probability as $\sup_{0\leq i\leq n-1}(t_{i+1}-t_i)\rightarrow 0$, and in this case the limit will be denoted by $\int_0^T y_t\circ dw_t.$
\end{definition}
\begin{Remark}\label{semimartingale}

%2.Suppose $(w_t)_{t\in[0,T]}$ is an $m$-dimensional Brownian motion. If $\textbf{y}=(y_t)_{t \in[0,T]}$ is an $n$-dimensional continuous semimartingale of the form
% $$y_t=y_0+\int_0^t \upsilon_sds + \sum_{i=1}^{m}\int_0^t  \zeta^{i}_s dw^{i}_s,$$
%where $(\upsilon_t)_{t \in[0,T]}$ and $(\zeta_t)_{t \in[0,T]}$ are adapted processes taking value in $\mathbb{R}^{n}$ and $\mathbb{R}^{n\times m}$ such that $\int_{0}^{T}|\upsilon_s|ds<\infty$ and $\int_{0}^{T}\zeta_s^2ds<\infty$ a.s..
%Then the Stratonovich integral of $\textbf{y}$ exits on any interval $[0,t]$, and
%$$\int_0^t y_s\circ dw_s^i=\int_0^t y_sdw^i_s+\frac{1}{2}\int_0^t  \textbf{1}_{m}^{\top}\zeta_s^i ds \cdot \textbf{1}_{m},$$
%where $\textbf{1}_{m}$ is the $m$-dimensional all-ones vector.

 We can translate an Ito integral to a Stratonovich integral. If $\textbf{y}=(y_t)_{t \in[0,T]}$ is a continuous semimartingale of the form
 $$y_t=y_0+\int_0^t \upsilon_s\,ds +\int_0^t  \zeta_s\,dw_s,$$
where $(\upsilon_t)_{t \in[0,T]}$ and $(\zeta_t)_{t \in[0,T]}$ are adapted processes taking value in $\mathbb{R}^{n}$ and $\mathbb{R}^{n\times m}$ such that\\
 $\int_{0}^{T}\parallel \upsilon_s\parallel ds<\infty$ and $\int_{0}^{T}\parallel\zeta_s\parallel^2ds<\infty$ a.s..
Then   $\textbf{y}$ is Stratonovich integrable on any interval $[0,t]$, and
\begin{align}
\int_0^t y_s\circ dw_s=\int_0^t y_s\,dw_s+\langle y,w\rangle_t=\int_0^t y_s\,dw_s+\frac{1}{2}\int_0^t \zeta_s\, ds ,\label{Stratonovich_equality}
\end{align}
where $\langle y,w \rangle_t$ denotes the joint quadrature variation of the semimartingale $\textbf{y}$ and the Brownian motion $\textbf{w}$.
Definition \ref{Stratonovich_integral} and the equality (\ref{Stratonovich_equality}) can be naturally extended to the vector case.% when $\{w_t\}$ is $m$-dimensional  and $\textbf{y}$ is  $n\times m$
\end{Remark}

%\begin{theorem}
%Consider a process of the form $$x_t=x_0+\int_0^t K_sds + \int_0^t  \Upsilon_sdw_s$$
%Then we have
%\begin{align*}
%F(t,x_t)=&F(0,x_0)+ \int_0^t F_s(s,x_s)ds+\int_0^t F_x(s,x_s)dx_s \\
%&+\frac{1}{2} \int_0^t F_{xx}(s,x_s)D_s^2ds+ \int_0^t F_xx(s,x_s)(D^{-}x)_s\Upsilon_sds.
%\end{align*}
%\end{theorem}
Then we present the Ito formula for  Stratonovich integral in Proposition \ref{Theorem7:Ito_Stratonovich} (see, e.g., Section 3.2.3 of \cite{nualart2006malliavin}).
\begin{Proposition}[Theorem 3.2.6 in \cite{nualart2006malliavin}]\label{Theorem7:Ito_Stratonovich}
Let $\textbf{w}=(w_t^1,\cdots,w_t^m)_{t\in[0,T]}$ be an $m$-dimensional Brownian motion. Suppose that $y_0\in \mathbb{D}^{1,2}$, $\upsilon_s \in \mathbb{L}^{1,2}$, and $\zeta^{i}\in \mathbb{L}_S^{2,4}$, $i=1,\cdots,m$. Consider a process $\textbf{y}=(y_t)_{t \in[0,T]}$ of the form
$$y_t=y_0+\int_0^t \upsilon_s\, ds + \sum_{i=1}^{m} \int_0^t  \zeta^{i}_s\circ dw^{i}_s, ~0\leq t\leq T.$$

Assume that $(y_t)_{0 \leq t\leq T}$ has continuous paths. Let $F:\mathbb{R}^{n}\rightarrow \mathbb{R}$ be a twice continuously differentiable function. Then we have
\begin{align}
F(y_t)=F(y_0)+ \int_0^t F^{\top}_y(y_s)\upsilon_s\, ds + \sum_{i=1}^{m}\int_0^t \left[F_y(y_s)^{\top}\zeta^i_s\right]\circ dw^{i}_s,  ~0\leq t\leq T,  \label{Stratonovich_ito_formula}
\end{align}
where $F_y(\cdot)$ denotes the gradient of $F$ w.r.t. $y$.
\end{Proposition}

Proposition \ref{Theorem7:Ito_Stratonovich} basically  says that the Stratonovich integral obeys the ordinary chain rule.
 %We can translate (\ref{Stratonovich_ito_formula}) into Ito language based on Remark \ref{semimartingale} to obtain the well-known Ito differential rule, when $\textbf{y}$ is a countinuous semimartingale.

Based on the definition of Stratonovich integral and Remark \ref{semimartingale}, we generalize the SDE (\ref{state_SDE}) to the  Stratonovich sense (referred to as S-SDE)  assuming that $b$ is bounded and $C^1$ in $(x,u)$; $\sigma$ is bounded and $C^2$ in  $x$. Then (\ref{state_SDE}) is equivalent to
\begin{align}
x_t=x+\int_0^t\bar{b}(t,x_t,u_t)dt+ \sum_{i=1}^{m} \int_0^t \sigma^i(t,x_t)\circ dw^i_t,~~0\leq t \leq T, \label{SDE_Anticipative1}
\end{align}
where  $\sigma^{i}:[0,T]\times\mathbb{R}^{n}\rightarrow\mathbb{R}^{n}$ is the $i$-th column of $\sigma$, $i=1,\cdots,m,$ and $\bar{b}(t,x,u)=b(t,x,u)- \frac{1}{2}\sum_{i=1}^{m}\sigma^i_x\sigma^i(t,x)$. Here $\sigma^i_x\sigma^i(t,x)$ denotes an $n\times 1$ vector with $\sum_{j=1}^{n}\frac{\partial \sigma^{ki}}{\partial x_j}(t,x)\sigma^{ji}(t,x)$ being its $k$-th entry and $\sigma^{ki}(\cdot)$ is the $k$-th component of $\sigma^i(\cdot)$. Since  the stochastic integral in (\ref{SDE_Anticipative1}) is in  the Stratonovich sense, S-SDE (\ref{SDE_Anticipative1}) adopts its solution in the space of $\mathcal{B}([0,T])\times\mathcal{F}$-measurable processes, which may not be adapted to the filtration generated by the Brownian motion. Therefore, we are  allowed to consider anticipative policies $\textbf{u}\in\mathcal{U}(0)$ in (\ref{SDE_Anticipative1}).  %Given $x_0=x$, equation (\ref{SDE_Anticipative1}) can be abbreviated as the equivalent Stratonovich stochastic

Finally, we need to ensure the existence of a solution to S-SDE (\ref{SDE_Anticipative1}) if the control strategy  $\textbf{u}\in \mathcal{U}(0)$ is anticipative. Following \cite{davis1992deterministic},\cite{Ocone1989}, we  have a representation of such a solution using the decomposition technique:
\begin{align}
x_t=\xi_t(\eta_t), \label{SDE_Anticipative_solution}
\end{align}
where $\{\xi_t(x)\}_{t\in[0,T]}$ denotes the stochastic flow defined by the adapted equation:
\begin{align}
d\xi_t&=\sum_{i=1}^{m} \sigma^{i}(t,\xi_t)\circ dw^i_t, \notag\\
&=\frac{1}{2}\sum_{i=1}^{m}\sigma^i_x\sigma^i(t,\xi_t)dt+ \sigma(t,\xi_t) dw_t,  ~~\xi_0=x, \label{SDE_decompostion1}
\end{align}
and $(\eta_t)_{t\in[0,T]}$ solves an ordinary differential equation:
\begin{equation}
\frac{d\eta_t}{dt}=\left(\frac{\partial \xi_t}{\partial x}\right)^{-1}(\eta_t)\bar{b}\left(t,\xi_t(\eta_t),u_t\right),  ~~\eta_0=x, \label{SDE_decompostion2}
\end{equation}
where $\frac{\partial \xi_t}{\partial x}$ denotes the $n\times n$ Jacobian matrix of $ \xi_t$ with respect to $x$. Under some technical conditions (see Section 1 of \cite{davis1992deterministic}), the solution (\ref{SDE_Anticipative_solution}) is defined almost surely: observe that $\xi_t$  does not depend on the control $u_t$, i.e., it is the solution to a regular SDE in the Ito sense; $\eta_t$ is not defined by a stochastic integral so it is the solution to an ordinary differential equation parameterized by $\textbf{w}$ (note that $\frac{\partial \xi_t}{\partial x}$ is well-defined a.s. for $(t,x)\in [0,T]\times \mathbb{R}^{n}$, because $\xi_t(x)$  is flow of diffeomorphisms a.s..). Hence, $x_t=\xi_t(\eta_t)$ is well-defined regardless of the adaptiveness of $\textbf{u}=(u_t)_{0\leq t \leq T}$. To check that $x_t=\xi_t(\eta_t)$ satisfies (\ref{SDE_Anticipative1}), we need to employ a generalized Ito formula of (\ref{Stratonovich_ito_formula}) for Stratonovich integral (see Theorem 4.1 in \cite{Ocone1989}).

%We come back to our problem.  Let $y_t=(t,x_t)$ and  let $B_i(y_t)=(0, \sigma_i(t,x_t)^{\top})^{\top}$ and $\bar{A}(y_t,u_t)=(1,(b(t,x_t,u_t)-\frac{1}{2}\sum_{i}^{m}\sigma_x^i\sigma^i(t,x_t))^{\top})^{\top}$, where $\sigma_i$ is the $i$-th column of $\sigma$.
\vspace{-3mm}
\subsection{Value Function-Based Penalty}\label{Section:Appendix:value_penalty}
The tools we have introduced in the last subsection, especially the  Ito formula for  Stratonovich integral, enable us to show the value function-based optimal penalty for the controlled Markov diffusions that developed in Theorem\ref{Theorem5:StrongDual_IdealPenalty}. % under the following assumption.

\begin{IEEEproof}[Proof of Theorem\ref{Theorem5:StrongDual_IdealPenalty}]  Suppose $\textbf{u}\in \mathcal{U}_{\mathbb{F}}(0)$ and let $y_t=V^{\top}_x(t,x_t)\sigma^i(t,x_t)$ in Remark \ref{semimartingale}  for $i=1,\cdots,m$. We can immediately obtain
\begin{align*}
h^{*}_v(\textbf{u},\textbf{w})=\sum_{i=1}^m\int_0^TV^{\top}_x(t,x_t)\sigma^i(t,x_t)\,dw_t^i=\int_0^TV^{\top}_x(t,x_t)\sigma(t,x_t)\,dw_t.
\end{align*}
%$h^{*}_v(\textbf{u},\textbf{w})=\sum_{i=1}^{m}\int_0^TV_x(t,x_t)^{\top}\sigma^i(t,x_t)dw^i_t=\int_0^TV_x(t,x_t)^{\top}\sigma(t,x_t)dw_t.$
Note that $V_x$ and $\sigma$ both satisfy a polynomial growth, since  $V(t,x)\in C^{1,2}(Q)\cap C_{p}(\bar{Q})$. Then  we have %and $\mathbb{E}[\int_{0}^{T}||u_t||^2dt]<\infty$,

$$\mathbb{E}_{0,x}\left[\parallel\int_{0}^{T}V^{\top}_{x}(t,x_t)\sigma(t,x_t) \parallel^2 dt\right]<\infty,$$
and therefore,  $\mathbb{E}_{0,x}[h^{*}_v(\textbf{u},\textbf{w})]=0$ when $\textbf{u}\in \mathcal{U}_{\mathbb{F}}(0)$. Hence, $h^{*}_v(\textbf{u},\textbf{w})\in  \mathcal{M}_{\mathbb{F}}(0)$.
We then show the strong duality
\begin{align}
V(0,x)=\mathbb{E}_{0,x}\left[\sup_{\textbf{u}\in \mathcal{U}(0)}\left\{\Lambda(x_{T})+\int_{0}^{T}g(t,x_t,u_t)dt-h^{*}_v(\textbf{u},\textbf{w})\right\}\right]. \label{strong_duality_proof}
\end{align}
According to the weak duality (i.e., Proposition \ref{Prop:weak_duality}),
\begin{align}
V(0,x)\leq \mathbb{E}_{0,x}\left[\sup_{\textbf{u}\in \mathcal{U}(0)}\left\{\Lambda(x_{T})+\int_{0}^{T}g(t,x_t,u_t)dt-h^{*}_v(\textbf{u},\textbf{w})\right\}\right]. \label{weak_duality1}
\end{align}
Next we prove the reverse inequality. Note that with $x_0=x$,
\begin{align*}
&\Lambda(x_{T})+\int_{0}^{T}g(t,x_t,u_t)dt-h^{\ast}_{v}(\textbf{u},\textbf{w})\\
=&~V(0,x)+\int_{0}^{T} \left[V_t(t,x_t)+V^{\top}_{x}(t,x_t)\bar{b}(t,x_t,u_t)\right]dt\\
&+\sum_{i=1}^{m}\int_{0}^{T}\left[V^{\top}_{x}(t,x_t)\sigma^i(t,x_t)\right]\circ dw^i_t-h^{\ast}_{v}(\textbf{u},\textbf{w})\\
=&~V(0,x)+\int_{0}^{T}\left[g(t,x_t,u_t)+A^{u_t}V(t,x_t)\right]dt,
\end{align*}
where the first equality is obtained by applying Ito formula for  Stratonovich integral (i.e., Proposition \ref{Theorem7:Ito_Stratonovich}) on $V(t,x)$ with $V(T,x_T)=\Lambda(x_T)$:
\begin{align*}
V(T,x_T)=&~V(0,x_0)+\int_{0}^{T} \left[V_t(t,x_t)+V^{\top}_{x}(t,x_t)\bar{b}(t,x_t,u_t)\right]dt\\
&+\sum_{i=1}^{m}\int_{0}^{T}\left[V^{\top}_{x}(t,x_t)\sigma^i(t,x_t)\right]\circ dw^i_t. %\label{Ito_rule}
\end{align*}
Since we assume the value function satisfies all the assumptions in Theorem \ref{Theorem2:Verification}(b), there exists an optimal control $\textbf{u}^{*}=(u^{*}_t)_{t\in[0,T]}$ with $u^{*}_t=u^{*}(t,x_t)$  and it satisfies
$$g(t,x,u^{\ast}(t,x))+A^{u^{*}(t,x)}V(t,x)=\max_{\textbf{u}\in \mathcal{U} }\left\{g(t,x,u)+A^{u}V(t,x)\right\}=0,$$
then we have
\begin{align}
&\sup_{\textbf{u}\in \mathcal{U}(0)}\left\{\Lambda(x_{T})+\int_{0}^{T}g(t,x_t,u_t)dt-h^{\ast}_{v}(\textbf{u},\textbf{w})\right\}\notag\\
=&\sup_{\textbf{u}\in \mathcal{U}(0)}\left\{V(0,x)+\int_{0}^{T}\bigg[g(t,x_t,u_t)+A^{u_t}V(t,x_t)\bigg]dt\right\}\notag\\
\leq&V(0,x)+ \int_{0}^{T} \sup_{u\in \mathcal{U}}\bigg\{g(t,x_t,u)+A^{u}V(t,x_t)\bigg\}dt \label{inequality_1}\\
=&V(0,x)+ \int_{0}^{T} \left[g(t,x^{*}_t,u^{*}_t)+A^{u_t^{*}}V(t,x_{t}^{*})\right]dt \notag\\
=&V(0,x). \label{inequality_2}
\end{align}
Taking the conditional expectation on both sides, we have
 $$V(0,x)\geq \mathbb{E}_{0,x}\left[\sup_{\textbf{u}\in \mathcal{U}(0)}\left\{\Lambda(x_{T})+\int_{0}^{T}g(t,x_t,u_t)dt-h^{\ast}_{v}(\textbf{u},\textbf{w})\right\}\right].$$
Together with the weak duality (\ref{weak_duality1}) , we reach the equality (\ref{strong_duality_proof}).

Due to the fact of the  equality (\ref{strong_duality_proof}) (that is in expectation sense) and the pathwise inequality (\ref{inequality_2}), we find that the only inequality (\ref{inequality_1}) should be an equality in almost sure sense. So the equality (\ref{almost_sure}) holds in almost sure sense. To achieve the equality in (\ref{inequality_1}), the optimal control $\textbf{u}^{*}$ should be applied, which implies the equality (\ref{almost_sure_2}).

%Due to the fact of the  equality (\ref{strong_duality}) (that is in expectation sense) and the pathwise inequality (\ref{inequality_2}), we find that the only inequality (\ref{inequality_1}) should be an equality in almost sure sense. So the equality (\ref{almost_sure}) holds in almost sure sense. To achieve the equality in (\ref{inequality_1}), the optimal control $\textbf{u}^{*}$ should be applied, which implies the equality (\ref{almost_sure_2}).
\end{IEEEproof}

%By imposing the value function-based optimal penalty the objective value of the dual problem is equal to $V(0,x)$ not only in the expectation sense, but also in the almost sure sense. Therefore, we can view the dual approach as a variance reduction technique. In particular,  $h^*_v$ plays the role of control variates.  As another obvious fact, $h^{*}_v(\textbf{u},\textbf{w})$ evaluated at $\textbf{u}=\textbf{u}^*$ is equal to $h^{*}(\textbf{u}^*,\textbf{w})$ in Proposition \ref{Prop:optimal_penalty}.

 %maximizing the inner optimization problem with respect to the anticipative control $\textbf{u}$ becomes an maximization problem over the penalized reward.
 %of this proposition is of the same intrisic of , and also relies on the Doob -Meyer Decompostion of the process $V(t,x_t)$ (that is in fact a supermartingale), which can be immediatedly derived by applying the Ito lemma on $V(t,x_t)$. $h^*_t$ in () is of the same form of  ,the difference that the control $u$  here is  $\max_{t\in[0,T]}$.  involved.  The reason is that the only control in the optimal stopping problem is to continue or stop the process but does not affect the evolution the process. The control in the dual approach is represented by $\max_{t\in[0,T]}$.

\vspace{-2mm}
 \subsection{Optimal Stopping under Diffusion Processes and Its Dual Representation}\label{Section:Appendix:optimal_stopping}
References \cite{Belomestny:2009,wang2010fast} use the martingale duality approach to compute upper bounds on the prices of American options, which is a typical optimal stopping problem. By viewing the martingale-based dual approach as a case of the perfect information relaxation, \cite{Belomestny:2009,wang2010fast} both explored the structure of the ``optimal penalty'' to the dual of the optimal stopping problem under the \emph{diffusion} process. We briefly review these results that parallel Theorem\ref{Theorem5:StrongDual_IdealPenalty} for controlled diffusions.

Suppose an uncontrolled diffusion $(x_t)_{t\in[0,T]}$ follows the SDE
\begin{align*}
dx_{t}&=b(t,x_t)dt+\sigma(t,x_t)dw_t,~~ 0 \leq t \leq T.
\end{align*}

We still use $\mathbb{F}$ to denote the natural filtration generated by the Brownian motion $(w_t)_{t\in[0,T]}.$ The primal representation of the optimal stopping problem is
\begin{align}
V(t,x)&=\sup_{\tau\in \mathcal{J}_t}\mathbb{E}_{t,x}[g(\tau, x_{\tau})], \label{primal_optimal_stopping_penalty}
\end{align}
where $g:\bar{Q}\rightarrow\mathbb{R}$  is a reward function,  and $\mathcal{J}_t$ is the set of $\mathbb{F}$-stopping times taking value in $[t,T]$.  Suppose that $V(t,x)$ is uniformly bounded and  is sufficiently smooth to apply Ito formula,  we have the following dual representation of the optimal stopping problem.
\begin{Proposition}[Theorem 1 and Theorem 2 in \cite{wang2010fast} ]
Let $\mathcal{H}_{\mathbb{F}}$ represent the space of $\mathbb{F}$-martingales $\{h_t\}_{t\in[0,T]}$ with $h_0$ = 0 and $\sup_{t\in[0,T]}\mathbb{E}[|h_t|]<\infty$. Then
\begin{equation}
V(0,x)=\min_{h\in \mathcal{H}_{\mathbb{F}}}\mathbb{E}_{0,x}\big[\max_{t\in[0,T]}\{g(t, x_{t})-h_t\}\big], \label{dual_optimal_stopping}
\end{equation}
In particular, the optimal martingale $\{h^*_t\}_{t\in[0,T]}$ that achieves the minimum in (\ref{dual_optimal_stopping}) is of the form
\begin{equation}
h^*_t=\int_{0}^{t}V_x(s,x_s)^{\top}\sigma(s,x_s)dw_s. \label{dual_optimal_stopping_penalty}
\end{equation}
\end{Proposition}
Noting that the maximization problem inside the expectation term (\ref{dual_optimal_stopping}) is the ``inner optimization problem'' in the dual representation of the optimal stopping problem,  since the only control in the primal (\ref{primal_optimal_stopping_penalty}) is to choose ``continue'' or ``stop''  the process. The strong duality result (\ref{dual_optimal_stopping}) holds for general Markov processes, which relies on the the Doob-Meyer decomposition of the process $\{V(t,x_t)\}_{t\in[0,T]}$; however,  the form of the optimal martingale (or penalty) $h^*$ in (\ref{dual_optimal_stopping_penalty}) is true only under the diffusion process. The form of $h^*$ exposes its connection with the value function-based penalty presented in Theorem\ref{Theorem5:StrongDual_IdealPenalty}.

\end{document}